\documentclass [leqno] {amsart}
  
\usepackage {amssymb}
\usepackage {amsthm}
\usepackage {amscd}

\usepackage[all]{xy}
\usepackage {hyperref} 

\begin{document}
\bibliographystyle{alpha}
\theoremstyle{plain}
\newtheorem{proposition}[subsubsection]{Proposition}
\newtheorem{lemma}[subsubsection]{Lemma}
\newtheorem{corollary}[subsubsection]{Corollary}
\newtheorem{thm}[subsubsection]{Theorem}
\newtheorem{introthm}{Theorem}
\newtheorem*{thm*}{Theorem}
\newtheorem{conjecture}[subsubsection]{Conjecture}
\newtheorem{question}[subsubsection]{Question}
\newtheorem{fails}[subsubsection]{Fails}

\theoremstyle{definition}
\newtheorem{definition}[subsubsection]{Definition}
\newtheorem{notation}[subsubsection]{Notation}
\newtheorem{condition}[subsubsection]{Condition}
\newtheorem{example}[subsubsection]{Example}
\newtheorem{claim}[subsubsection]{Claim}

\theoremstyle{remark}
\newtheorem{remark}[subsubsection]{Remark}

\numberwithin{equation}{subsection}

%Matheumgebungen
\newcommand{\eq}[2]{\begin{equation}\label{#1}#2 \end{equation}}
\newcommand{\ml}[2]{\begin{multline}\label{#1}#2 \end{multline}}
\newcommand{\mlnl}[1]{\begin{multline*}#1 \end{multline*}}
\newcommand{\ga}[2]{\begin{gather}\label{#1}#2 \end{gather}}
\newcommand{\mat}[1]{\left(\begin{smallmatrix}#1\end{smallmatrix}\right)}

%xypic
\newcommand{\arir}{\ar@{^{(}->}}
\newcommand{\aril}{\ar@{_{(}->}}
\newcommand{\are}{\ar@{>>}}

% Pfeile xr , xl 
\newcommand{\xr}[1] {\xrightarrow{#1}}
\newcommand{\xl}[1] {\xleftarrow{#1}}
\newcommand{\lra}{\longrightarrow}
\newcommand{\inj}{\hookrightarrow}

% mathfrac, mathcal 
\newcommand{\mf}[1]{\mathfrak{#1}}
\newcommand{\mc}[1]{\mathcal{#1}}

% rm - Abkuerzungen 
\newcommand{\CH}{{\rm CH}}
\newcommand{\Gr}{{\rm Gr}}
\newcommand{\codim}{{\rm codim}}
\newcommand{\cd}{{\rm cd}}
\newcommand{\Spec} {{\rm Spec}}
\newcommand{\supp} {{\rm supp}}
\newcommand{\Hom} {{\rm Hom}}
\newcommand{\End} {{\rm End}}
\newcommand{\id}{{\rm id}}
\newcommand{\Aut}{{\rm Aut}}
\newcommand{\sHom}{{\rm \mathcal{H}om}}
\newcommand{\Tr}{{\rm Tr}}
\newcommand{\lf}{[\![}
\newcommand{\rf}{]\!]}

% Abkuerzungen zu Bourbaki-Notation

% Zahlen
\renewcommand{\P} {\mathbb{P}}
\newcommand{\Z} {\mathbb{Z}}
\newcommand{\Q} {\mathbb{Q}}
\newcommand{\C} {\mathbb{C}}
\newcommand{\F} {\mathbb{F}}

%sonst
\newcommand{\OO}{\mathcal{O}}

\newcommand{\qc}{q-\mathrm{crys}}
\newcommand{\dd}{[p]_q}
\newcommand{\wotimes}{\widehat{\otimes}}
\newcommand{\tR}{\tilde{R}}
\newcommand{\qdiff}{q{\rm -HPDiff}}

%test
 
\title{$q$-crystals and $q$-connections}
\author{Andre Chatzistamatiou}

%\thanks{This work has been supported by my Heisenberg fellowship}

\begin{abstract}
We study how the category of $q$-connections depends on the choice of coordinates. We exploit Bhatt's and Scholze's $q$-crystalline site, which is based on a coordinate free formulation of $q$-PD structures, in order to relate $q$-crystals and $q$-connections in the $p$-adic setting. This yields a natural equivalence between the the categories of $q$-connections for different choices of coordinates in the $p$-adic setting. The equivalence can be described explicitly in terms of differential operators. 

In order to obtain a global equivalence, we patch the $p$-adic differential operators to create a global one. The process is not entirely formal, and we are only able to obtain a global equivalence after inverting $2$.   
\end{abstract}

\maketitle

%\tableofcontents

\section*{Introduction}

The subject of this paper are Scholze's conjectures about $q$-connections \cite{S}. In their simplest form $q$-connections are modules, satisfying a completeness condition, over a certain subalgebra of the ring of differential operators. 

For example, let $R_q=\Z[x]\lf q-1 \rf$, and let $R_q\{\nabla_{x,q}\}$ be the $(q-1)$-completion of the non commutative $R_q$-algebra generated by $\nabla_{x,q}$ and satisfying $\nabla_{x,q}\cdot x - q\cdot x\cdot \nabla_{x,q} = 1$. This is a $q$-deformation of the Weyl algebra $\Z[x]\{\nabla_{x}\}$. It is an arithmetic phenomenon that this deformation depends on the coordinate $x$. Replacing $\Z$ by $\Q$ would yield the trivial deformation of the Weyl algebra.

The category of $(q-1)$-(derived) complete $R_q\{\nabla_{x,q}\}$-modules is the category of $q$-connections over $\Z[x]$ with respect to the coordinate $x$. According to the conjectures of \cite{S}, this category does not depend on the choice of the coordinate $x$. For example, there should be a natural equivalence between the category of $(q-1)$-(derived) complete $R_q\{\nabla_{x,q}\}$-modules and $(q-1)$-(derived) complete $R_q\{\nabla_{x-1,q}\}$-modules, where $\nabla_{x-1,q} \cdot (x-1) - q\cdot (x-1)\cdot \nabla_{x-1,q} = 1$.

Our paper exploits the notion of $q$-crystals introduced in \cite{BS}. Our first equivalence is between $q$-crystals and $q$-HPD stratifications (Theorem \ref{theorem-equivalence}). Both notions specialize to the respective classical notions in crystalline theory for $q=1$, and the proof of the equivalence follows from the crystalline formalism. The main observation, which allows us to run the crystalline formalism, is the existence of local retractions in the $q$-crystalline site (Proposition \ref{proposition-retractions}). 

By definition, $q$-HPD stratifications depend on the choice of a ``lifting". Independence of this choice is established via the identification with $q$-crystals.

The second equivalence is between $(p,q-1)$-completely flat $q$-HPD stratifications and $(p,q-1)$-completely flat quasi-nilpotent $q$-connections (Proposition \ref{proposition-equivalence-stratification-connections}). Indirectly, it shows the independence of the choice of coordinates for $(p,q-1)$-completely flat quasi-nilpotent $q$-connections. By using differential operators, we can make this independence explicit and extend it to all $q$-connections (Proposition \ref{proposition-morphism-qdiff}). All this takes place in the $p$-adic setting and works for all primes.

In order to briefly explain how differential operators are used to transform from coordinates $\underline{x}$ to coordinates $\underline{y}$, let us suppose that $q$-connections for $\underline{x}$ (and $\underline{y}$) are $\mc{A}_{p,\underline{x}}$-modules (and $\mc{A}_{p,\underline{y}}$-modules, respectively). Moreover, $\mc{A}_{p,\underline{x}}$ and $\mc{A}_{p,\underline{y}}$ are subalgebras of $\mc{D}_p$, (a completion of) the ring of differential operators. We can show that there is a distiguished free left $\mc{A}_{p,\underline{x}}$-module and free right $\mc{A}_{p,\underline{y}}$-module $\mc{M}_{p,\underline{x},\underline{y}}$ contained in $\mc{D}_p$. A choice of a common generator $D\in \mc{M}_{p,\underline{x},\underline{y}}$ yields an isomorphism $\psi_D: \mc{A}_{p,\underline{x}} \xr{} \mc{A}_{p,\underline{y}}$ defined by $a\cdot D = D\cdot \psi_D(a)$. This isomorphism is well-defined up to inner automorphisms, and induces an equivalence between $\mc{A}_{p,\underline{x}}$-modules and $\mc{A}_{p,\underline{y}}$-modules up to natural transformation.

To obtain a global result, we use a natural condition to identify the differential operators $\mc{M}^{{\rm global}}_{p,\underline{x},\underline{y}} \subset \mc{M}_{p,\underline{x},\underline{y}}$ that are the image of global ones. This approach does not work for $p=2$, because $\mc{M}^{{\rm global}}_{2,\underline{x},\underline{y}}=\emptyset$ in general.

Finally, after inverting $2$, we patch $\{\mc{M}^{{\rm global}}_{p,\underline{x},\underline{y}}\}_{p>2}$ together to obtain a free left $\mc{A}_{\underline{x}}$-module and free right $\mc{A}_{\underline{y}}$-module $\mc{M}_{\underline{x},\underline{y}}$ contained in the ring of global differential operators $\mc{D}$.

\section{$q$-crystals}

This section heavily relies on \cite[\textsection16]{BS}. Bhatt and Scholze introduce the notion of a $q$-divided power algebra \cite[Definition~16.2]{BS}, a $q$-divided power envelope \cite[Lemma~16.10]{BS} and the $q$-crystalline site \cite[Definition~16.12]{BS}. These are major innovations which allows them to glue a $q$-de Rham complex, prove some of the conjectures in \cite{S}, and to relate their previous work on $A\Omega$-complexes to prismatic cohomology.  

\subsection{The $q$-crystalline site}

We fix a $q$-PD pair $(D,I)$ and a $p$-completely smooth and $p$-complete $D/I$-algebra $R$. 

\begin{definition}[The $q$-crystalline site] (See \cite[Definition~16.12]{BS})
The {\em $q$-crystalline site of $R$ relative to $D$}, denoted $(R/D)_{\qc}$, is the category of $q$-PD thickenings %of $p$-completed localizations 
of $R$ relative to $D$, i.e. the category of $q$-PD pairs $(E,J)$ over $(D,I)$ equipped with an isomorphism $R\xr{} E/J$ of $D/I$-algebras.
%$R[f^{-1}]^{\wedge} \xr{} E/J$ of $D/I$-algebras 
%for some $f\in R$.  
\end{definition}

An important property of the classical crystalline site and the infinitesimal site for a smooth scheme is the existence of local retractions (see \cite[\textsection4.2]{G}). This property is used to relate crystals to more explicit objects like stratifications and connections.

\begin{proposition}[Existence of retractions]\label{proposition-retractions}
Suppose there is $(\tilde{R}, I\tilde{R})$ in $(R/D)_{\qc}$ such that 
%$\tilde{R}/I\tilde{R}=R$ and 
a $(p, \dd)$-\'etale map $\psi:D[x_1,\dots,x_n]^{\wedge}\xr{} \tilde{R}$ exists. 

For every $(S,J)$ in $(R/D)_{\qc}$, the coproduct of $(S,J)$ and $(\tilde{R}, I\tilde{R})$ in   $(R/D)_{\qc}$ exists. We denote it by $(\tilde{S}, \tilde{J}) = (S,J) \otimes_{\qc} (\tilde{R}, I\tilde{R})$. The natural morphism  
$$
(S,J) \xr{f} (\tilde{S}, \tilde{J}) 
$$ 
is $(p, \dd)$-faithfully flat. Moreover, there is a section $\tilde{S}\xr{} S$ in the category of $S$-modules.
\begin{proof}
We note that $A:=D[x_1,\dots,x_n]^{\wedge}$, $\tilde{R}$, and $S$ are classically $(p,\dd)$-complete \cite[Lemma~3.7(1)]{BS}.
  
We can find a morphism of $D$-algebras $\tau': A\xr{} S$ lifting the map $A/I \xr{\psi} \tilde{R}/I = R \xr{} S/J$. As a first step, we will extend $\tau'$ to $\tilde{R}$. 

We can find a presentation $\tilde{R}=A[y_1,\dots,y_d]/(f_1,\dots,f_d)^{\wedge}$ such that ${\det}(\partial f_j/\partial y_i)$ is invertible. For $B=A[y_1,\dots,y_d]\lf z_1, \dots z_d \rf/(f_1-z_1,\dots,f_d-z_d)$, there is ring homomorphism $B\xr{} S$ extending $\tau':A\xr{} S$, mapping $z_i$ into $J$, and reducing to $R\xr{} S/J$. Indeed, for every $a\in J$, we have $a^p \in (p, \dd)$. After mapping all $y_i$ to some lifts, we map $z_i$ to the image of $f_i$.

Since $A[y_1,\dots,y_d]/(f_1,\dots,f_d)$ is \'etale over $A$, we can find a morphism of $A$-algebras $A[y_1,\dots,y_d]/(f_1,\dots,f_d) \xr{} B$ such that the composition with $B\xr{} B/(z_1,\dots,z_d)$ is the identity. The composition with $B\xr{} S$ yields the desired lift $\tau': \tilde{R} \xr{} S$ after completion.

We note that $\tau'$ is in general not compatible with the $\delta$-structures. However, there is a $\delta$-structure on $S\lf \epsilon_1, \dots, \epsilon_n\rf$ such that $\tau: \tilde{R} \xr{} S \lf \epsilon_1, \dots, \epsilon_n\rf$ defined by $x_i\mapsto \tau'(x_i)+\epsilon_i$ is a morphism of $\delta$-rings. In order to define this $\delta$-structure we need to solve $\delta(\tau(x_i))=\tau(\delta(x_i))$. Since
\begin{align*}
\delta(\tau(x_i)) &= \delta(\tau'(x_i)+\epsilon_i)=\delta(\tau'(x_i)) + \delta(\epsilon_i) + \sum_{i=1}^{p-1} p^{-1}\binom{p}{i} \tau'(x_i)^{i}\epsilon_i^{p-i},\\
\tau(\delta(x_i)) &= \sum_{I} \tau'(\delta(x_i)^{(I)}) \cdot \epsilon_1^{I_1}\cdots \epsilon_n^{I_n}, 
\end{align*}
where the second equation is the Taylor expansion and $\delta(x_i)^{(I)}=(\prod_{j=1}^n \frac{\partial_{x_j}^{I_j}}{I_j!})(\delta(x_i))$, we have to set
$$
\delta(\epsilon_i) = \sum_{I} \tau'(\delta(x_i)^{(I)}) \cdot \epsilon_1^{I_1}\cdots \epsilon_n^{I_n} - \delta(\tau'(x_i)) - \sum_{i=1}^{p-1} p^{-1}\binom{p}{i} \tau'(x_i)^{i}\epsilon_i^{p-i}.
$$
This extends to a $\delta$-structure on $S \lf \epsilon_1, \dots, \epsilon_n\rf$ with the desired properties. 

Finally, applying the $q$-PD envelope \cite[Lemma~16.10]{BS} construction over $(S,J)$ to the $(p, \dd)$-completely flat $\delta$-$S$-algebra $S \lf \epsilon_1, \dots, \epsilon_n\rf$ and the regular sequence $\epsilon_1,\dots, \epsilon_n$, we obtain $(\tilde{S}, \tilde{J})$. We denote by $f:(S,J)\xr{} (\tilde{S},\tilde{J})$ and $\tau: (\tilde{R}, I\tilde{R}) \xr{} (\tilde{S}, \tilde{J})$ the induced morphisms in $(R/D)_{\qc}$.

In order to show that $(\tilde{S}, \tilde{J})$ is the coproduct of $(S, J)$ and $(\tilde{R}, I\tilde{R})$, we can use the universal property of a $q$-PD envelope. For morphisms $g_{S}: S\xr{} T$ and $g_{\tilde{R}}: \tilde{R} \xr{} T$ in $(R/D)_{\qc}$, it suffices to show that there is a unique morphism of $\delta$-rings $g:S \lf \epsilon_1, \dots, \epsilon_n\rf \xr{} T$ such that $g\circ \tau=g_{\tilde{R}}$, $g\circ f = g_{S}$, and which maps the ideal $(\epsilon_1, \dots, \epsilon_n)$ into the $q$-PD ideal of $T$. Only the morphism of $S$-algebras induced by $\epsilon_i \mapsto g_{\tilde{R}}(x_i)-g_S(\tau'(x_i))$ satisfies these properties.   

We still need to show the existence of a section $\tilde{S}\xr{} S$ in the category of $S$-modules. 
Let $\gamma_p(a)=\frac{\phi(a)}{\dd}-\delta(a)$, we have $\gamma_p(\tilde{J})\subset \tilde{J}$. 
We set $\gamma_1={\rm id}_{\tilde{S}}$, $\gamma_{p^k}=\gamma_p\circ \gamma_{p^{k-1}}$, for all $k\geq 1$, and for any positive integer $i$ with $p$-adic expansion $i=\sum_{k=0}^{\infty} i_k p^k$, 
we set $\gamma_i(a) = \prod_{k=0}^{\infty} \gamma_{p^k}(a)^{i_k}$. This defines maps $\gamma_i: \tilde{J}\xr{} \tilde{J}$. 

For each $I=(i_1,\dots, i_n)\in \mathbb{Z}_{\geq 0}^n$, we set 
\begin{equation}\label{equation-gamma-I}
\Gamma_I:= \prod_{\substack{j=1\\i_j\neq 0}}^{n} \gamma_{i_j}(\epsilon_j),
\end{equation}
and claim that  
\begin{equation}\label{equation-lots-of-sections}
\left(\bigoplus_{I\in \mathbb{Z}_{\geq 0}^n} S \cdot \Gamma_I \right)^{\wedge} \xr{} \tilde{S},
\end{equation}
is an isomorphism of $S$-modules, where the left hand side is the derived $(p,\dd)$-completion. Indeed, $(p,\dd)$-completion is the same as $(p,q-1)$-completion, and by the derived Nakayama lemma it suffices to show the isomorphism after derived base change $\otimes^L_{S} S/(p,q-1)$. Both modules are $(p,\dd)$-completely flat, hence derived base change reduces to simple base change $\otimes_S S/(p,q-1)$. After reduction modulo $q-1$ the map already becomes an isomorphism. Indeed, $\tilde{S}/q-1$ is the $p$-completion of the pd-envelope of $S/(q-1)\lf \epsilon_1, \dots, \epsilon_n\rf$ along $(\epsilon_1, \dots, \epsilon_n)$, and 
$$
\Gamma_I \equiv \prod_j \frac{\epsilon_j^{i_j}}{[i_j!]_p} \mod q-1,
$$ 
where $[n]_p$ denotes the largest $p$-power dividing $n$.

By using \eqref{equation-lots-of-sections} and projecting to the summand corresponding to $I=0$, we obtain the desired section.
\end{proof}
\end{proposition}

\subsection{$q$-crystals}
For a $D$-algebra $E$, we have the abelian category of derived $(p,\dd)$-complete $E$-modules  ${\rm Mod}^{\wedge}_E$ at our disposal (see Appendix \ref{appendix-descent}). For a morphism of $D$-algebras $f:E\xr{} E'$, there is a right exact base change functor  
$$
f^*:{\rm Mod}^{\wedge}_E \xr{} {\rm Mod}^{\wedge}_{E'}, \quad M\mapsto E'\widehat{\otimes}_E M. 
$$

\begin{definition}[$q$-crystals]
%Let $\mathcal{C}$ be a full subcategory of $(R/D)_{\qc}$. 

A {\em $q$-crystal on $(R/D)_{\qc}$} is a derived $(p,\dd)$-complete $E$-module $M_E$ for each $(E,J)\in (R/D)_{\qc}$ together with isomorphisms 
$$
M_f:f^*M_E \xr{} M_{E'}, 
$$
for each $f: (E,J)\xr{} (E',J')$ in $(R/D)_{\qc}$ that satisfy 
$$
M_{g\circ f} =  M_g\circ g^*(M_f) %\left[ E''\widehat{\otimes}_E M_E \xr{E'' \widehat{\otimes} M_f} %E''\widehat{\otimes}_{E'} M_{E'} \xr{M_g} M_{E''} \right] 
$$
for all $f: (E,J)\xr{} (E',J')$ and $g: (E',J') \xr{} (E'', J'')$.

A {\em morphism of $q$-crystals} $u$ consists of a morphism of derived $(p,\dd)$-complete $E$-modules 
$$
u_E: M_E \xr{} M'_{E}
$$ 
for every $(E,J)$ in $(R/D)_{\qc}$ such that 
$$
u_{E'}\circ M_f = M'_f \circ  f^*(u_E)
$$
for all $f: (E,J)\xr{} (E',J')$.

We denote the category of $q$-crystals by $q{\rm -Cris}(R/D)$.

A $q$-crystal $M$ is called {\em $(p,\dd)$-completely flat} if $M_E$ is $(p,\dd)$-completely flat for every $(E,J)$ in $(R/D)_{\qc}$ (see Appendix \ref{appendix-completely-flat} for the definition of $(p,\dd)$-completely flat modules).
\end{definition}

\subsection{}\label{subsection-stratifiactions}
Let $\tilde{R}$ be as in Proposition \ref{proposition-retractions}. As a first step towards understanding $q{\rm -Cris}(R/D)$, we will outline which objects $q$-crystals induce on $\tilde{R}$. This will lead to the definition of {\em $q$-HPD stratifications}. The main result in this section will be an equivalence of categories between $q$-crystals and $q$-HPD stratifications if coordinates exist. This is analogous to the equivalence for the crystalline theory \cite[\textsection6]{BO}. 

Consider the diagram 
\begin{equation} \label{equation-delta-diagram}
\xymatrix{ \tilde{R} \ar@<1ex>[r]^-{\delta^1_0} \ar@<-1ex>[r]_-{\delta_1^1} &  \tilde{R}\lf \epsilon_1,\dots,\epsilon_n\rf  \ar@<3ex>[r]^-{\delta^2_0} \ar@<0ex>[r]^-{\delta^2_1} \ar@<-3ex>[r]^-{\delta^2_2} &  \tilde{R}\lf \epsilon_1,\dots,\epsilon_n, \tau_1, \dots, \tau_n \rf},
\end{equation}
where $\delta^1_0(x_i)=x_i+\epsilon_i$, $\delta^1_1(x_i)=x_i$, and 

\begin{align*}
&\delta^2_0(x_i) = x_i+\epsilon_i & &\delta^2_1(x_i) = x_i & &\delta^2_2(x_i) = x_i \\
&\delta^2_0(x_i+\epsilon_i) = x_i+\tau_i &  &\delta^2_1(x_i+\epsilon_i) = x_i + \tau_i & &\delta^2_2(x_i+\epsilon_i) = x_i +\epsilon_i. 
\end{align*}
We equip $\tilde{R}\lf \epsilon_1,\dots,\epsilon_n\rf$ and $\tilde{R}\lf \epsilon_1,\dots,\epsilon_n, \tau_1, \dots, \tau_n \rf$ with the $\delta$-structure that makes all maps in the diagram to maps of $\delta$-rings (see the proof of Proposition \ref{proposition-retractions} for how to make this work). 

Next, we apply the $q$-PD envelope \cite[Lemma~16.10]{BS} construction over $(\tilde{R},I\tilde{R})$ to the $(p, \dd)$-completely flat $\delta$-$\tilde{R}$-algebras $\tilde{R}\lf \epsilon_1,\dots,\epsilon_n\rf$ and $\tilde{R}\lf \epsilon_1,\dots,\epsilon_n, \tau_1, \dots, \tau_n \rf$ for the regular sequences $\epsilon_1,\dots, \epsilon_n$ and $\epsilon_1,\dots,\epsilon_n, \tau_1, \dots, \tau_n$, respectively. We denote the resulting objects in $(R/D)_{\qc}$ by $\tilde{R}^{(2)}$ and $\tilde{R}^{(3)}$. By using the universal property of a $q$-PD envelope we obtain an induced diagram 
\begin{equation}\label{equation-R-escalation}
\xymatrix{ \tilde{R} \ar@<1ex>[r]^-{\delta^1_0} \ar@<-1ex>[r]_-{\delta_1^1} &  \tilde{R}^{(2)}  \ar@<3ex>[r]^-{\delta^2_0} \ar@<0ex>[r]^-{\delta^2_1} \ar@<-3ex>[r]^-{\delta^2_2} &  \tilde{R}^{(3)}.}
\end{equation}
We have $\tilde{R}^{(2)}=\tilde{R}\otimes_{\qc} \tilde{R}$ and $\tilde{R}^{(3)}=\tilde{R}\otimes_{\qc} \tilde{R} \otimes_{\qc} \tilde{R}$. 

\begin{definition}
A  {\em $q$-HPD stratifications} on $(\tilde{R}, I\tilde{R})$ is a $(p, \dd)$-derived complete module $\tilde{R}$-module $M$ together with an isomorphism of $\tilde{R}^{(2)}$-modules 
$$
\epsilon: M\wotimes_{\tilde{R}, \delta^1_0} \tilde{R}^{(2)} \xr{} M\wotimes_{\tilde{R}, \delta^1_1} \tilde{R}^{(2)}
$$  
such that the following cocycle condition is satisfied:  
$$
\epsilon\wotimes_{\tilde{R}^{(2)}, \delta^2_2} \tilde{R}^{(3)} \circ \epsilon \wotimes_{\tilde{R}^{(2)}, \delta^2_0} \tilde{R}^{(3)} = \epsilon\wotimes_{\tilde{R}^{(2)}, \delta^2_1} \tilde{R}^{(3)}. 
$$
We form the category of $q$-HPD stratifications in the usual way. 

A $q$-HPD stratification $(M,\epsilon)$ is called {\em $(p,\dd)$-completely flat} if $M$ is a $(p,\dd)$-completely flat $\tilde{R}$-module (see Appendix \ref{appendix-completely-flat} for the definition of $(p,\dd)$-completely flat modules).
\end{definition}

\begin{remark}
It follows from the cocycle condition and the existence of an inverse for $\epsilon$ that  $ 
{\rm id}_M=\epsilon\wotimes_{\tilde{R}^{(2)}, m} \tilde{R}
$, where $m: \tilde{R}\otimes_{\qc} \tilde{R} \xr{} \tilde{R}$ is induced by the identity on each factor.
\end{remark}

We have a functor 
\begin{equation}\label{equation-crystals-to-connections}
q{\rm -Cris}(R/D) \xr{} \text{($q$-HPD stratifications on $(\tilde{R}, I\tilde{R})$)},
\end{equation}
defined by 
$$
M=(E\mapsto M_E, f\mapsto M_f) \mapsto (M_{\tilde{R}}, M_{\delta^1_1}^{-1} \circ M_{\delta^1_0}).
$$

\begin{thm}\label{theorem-equivalence}
Let $(\tilde{R}, I\tilde{R})$ be as in Proposition \ref{proposition-retractions}. The functor \eqref{equation-crystals-to-connections} from $q$-crystals to $q$-HPD stratifications on $\tilde{R}$ is an equivalence of categories.
\begin{proof}
For $(S,J)$ in $(R/D)_{\qc}$, Proposition \ref{proposition-retractions} yields morphisms 
\begin{equation} \label{equation-covering-S}
(S,J) \xr{f} (\tilde{S}, \tilde{J}) \xl{\tau} (\tilde{R}, I\tilde{R})
\end{equation}
in $(R/D)_{\qc}$. We define the diagram 
\begin{equation} \label{equation-descent-diagram}
\xymatrix{ \tilde{S} \ar@<1ex>[r]^-{\delta^1_0} \ar@<-1ex>[r]_-{\delta_1^1} &  \tilde{S}\wotimes_S \tilde{S}  \ar@<3ex>[r]^-{\delta^2_0} \ar@<0ex>[r]^-{\delta^2_1} \ar@<-3ex>[r]^-{\delta^2_2} &  \tilde{S}\wotimes_S \tilde{S} \wotimes_{S} \tilde{S}}
\end{equation}
in the usual way, that is, $\delta^1_0(s)=1\wotimes s$, $\delta^1_1(s)=s\wotimes 1$, and   
\begin{equation*}
\delta^2_0(s_1\wotimes s_2) = 1\wotimes s_1\wotimes s_2, \quad \delta^2_1(s_1\wotimes s_2) = s_1\wotimes 1\wotimes s_2, \quad\delta^2_2(s_1\wotimes s_2) = s_1\wotimes s_2\wotimes 1. 
\end{equation*}
%Note that $(\tilde{S}\wotimes_S \tilde{S}, \tilde{J}\wotimes_{S} \tilde{S} + \tilde{S}\wotimes_{S} \tilde{J})$ and $(\tilde{S}\wotimes_S \tilde{S} \wotimes_{S} \tilde{S}, \tilde{J}\wotimes_{S} \tilde{S} + \tilde{S}\wotimes_{S} \tilde{J} \wotimes_{S} \tilde{S} + \tilde{S} \wotimes_{S} \tilde{S} \wotimes_{S} \tilde{J})$ are objects in $(R/D)_{\qc}$.
Note that $\tilde{S}\wotimes_S \tilde{S}$ and  $\tilde{S}\wotimes_S \tilde{S} \wotimes_{S} \tilde{S}$ are objects in $(R/D)_{\qc}$. They are $q$-PD thickenings of $S/J$. 

Next, we extend $\tau:\tilde{R} \xr{} \tilde{S}$ to morphisms $\tau^{(2)}: \tilde{R}^{(2)}\xr {} \tilde{S}\wotimes_S \tilde{S}$ and $\tau^{(3)}: \tilde{R}^{(3)}\xr {} \tilde{S}\wotimes_S \tilde{S}\wotimes_S \tilde{S}$ in $(R/D)_{\qc}$. We want 
\begin{align*}
\tau^{(2)} \circ \delta_{1}^1 &= \delta_{1}^1 \circ \tau,  & \tau^{(2)}(\epsilon_i) &=   1\wotimes \tau(x_i) - \tau(x_i)\wotimes 1, & & \\
\tau^{(3)} \circ \delta_{2}^2 &= \delta_{2}^2 \circ \tau^{(2)}, &
\tau^{(3)}(\tau_i) &=   1\wotimes 1\wotimes \tau(x_i) - 1\wotimes \tau(x_i)\wotimes 1.
\end{align*}
This induces well-defined maps, because $\tilde{S}/\tilde{J}=S/J$ implies $1\wotimes \tau(x_i) - \tau(x_i)\wotimes 1 \in \tilde{J}\wotimes_{S} \tilde{S} + \tilde{S}\wotimes_{S} \tilde{J}$. 

By using $\tau, \tau^{(2)},$ and $\tau^{(3)}$, we get a morphism of diagrams $\text{\eqref{equation-R-escalation}}\xr{} \text{\eqref{equation-descent-diagram}}$. And this yields a functor from $q$-HPD stratifications to descent data
$$
\text{($q$-HPD stratifications on $(\tilde{R}, I\tilde{R})$)} \xr{} DD^{\wedge}_{\tilde{S}/S}
$$
(see Appendix \ref{appendix-descent}). Proposition \ref{proposition-retractions} guarantees a section $\tilde{S}\xr{} S$ in the category of $S$-modules. Proposition \ref{proposition-section-descent} implies ${\rm Mod}^{\wedge}_{S} \cong DD^{\wedge}_{\tilde{S}/S}$. We will denote by $(f,\tau)^*$ the resulting functor 
$$
\text{($q$-HPD stratifications on $(\tilde{R}, I\tilde{R})$)} \xr{} {\rm Mod}^{\wedge}_{S}.
$$

The next step is to define natural isomorphisms 
\begin{equation}\label{equation-natural-transformations-descent}
T_{(f_0,\tau_0)^*, (f_1,\tau_1)^*} : (f_0,\tau_0)^* \xr{\cong} (f_1,\tau_1)^*  
\end{equation}
for two different choices for the retractions in \eqref{equation-covering-S}. We will assume that both choices have a section  $\tilde{S_i}\xr{} S$, $i=0,1$. Taking $\tilde{S} = \tilde{S}_1\wotimes_S \tilde{S}_0$, $f = f_1\wotimes 1 = 1\wotimes f_0$, and identifying $DD^{\wedge}_{\tilde{S}_1/S}\cong  DD^{\wedge}_{\tilde{S}/S} \cong DD^{\wedge}_{\tilde{S}_0/S}$ we reduce to constructing
$$
(f,1\wotimes \tau_0)^* \xr{\cong} (f, \tau_1\wotimes 1)^*.  
$$  
To simplify the notation, we will simply write $\tau_1$ for $\tau_1\wotimes 1$, and 
similarly for $\tau_0$.

We can define $\mu:\tilde{R}^{(2)}\xr{} \tilde{S}$ in $(R/D)_{\qc}$ such that $\mu\circ \delta^1_{0}=\tau_0$ and $\mu\circ \delta^1_{1}=\tau_1$. This gives a natural isomorphism 
$$
\epsilon\wotimes_{\tilde{R},\mu} \tilde{S}: M\wotimes_{\tilde{R}, \tau_0} \tilde{S} \xr{} M\wotimes_{\tilde{R}, \tau_1} \tilde{S}
$$ 
for all $q$-HPD stratifications $(M,\epsilon)$. We claim that this induces an isomorphism of descent data. In other words, we have to prove the equality
\begin{equation}\label{equation-morphism-of-descent-data}
(\epsilon \wotimes_{\tilde{R}^{(2)}, \tau_1^{(2)}} \tilde{S} \wotimes_S \tilde{S}) \circ (\epsilon\wotimes_{\tilde{R}^{(2)},\delta^1_{0}\circ \mu} \tilde{S} \wotimes_S \tilde{S}) = (\epsilon\wotimes_{\tilde{R}^{(2)},\delta^1_{1}\circ \mu} \tilde{S} \wotimes_S \tilde{S}) \circ (\epsilon \wotimes_{\tilde{R}^{(2)}, \tau_0^{(2)}} \tilde{S} \wotimes_S \tilde{S}).
\end{equation}
We can define morphisms in $(R/D)_{\qc}$:
$$
\xymatrix{ \tilde{R}^{(3)} \ar@<1ex>[r]^-{\mu_0} \ar@<-1ex>[r]_-{\mu_1} &  \tilde{S}\wotimes_S \tilde{S}}
$$
such that $\mu_0\circ \delta_0^2 = \tau_0^{(2)}$, $\mu_{0}\circ \delta_2^2=\delta^1_1\circ \mu$, $\mu_1\circ \delta_2^2 = \tau_1^{(2)}$, $\mu_1\circ \delta^2_0=\delta^1_0\circ \mu$, and $\mu_0\circ \delta^2_1 = \mu_1\circ \delta^2_1$. After applying $\wotimes_{\tilde{R}^{(3)},\mu_0} \tilde{S}\wotimes_S \tilde{S}$ and $\wotimes_{\tilde{R}^{(3)},\mu_1} \tilde{S}\wotimes_S \tilde{S}$ to the cocycle condition, we obtain \eqref{equation-morphism-of-descent-data}. 

At this point we have constructed the natural isomorphisms \eqref{equation-natural-transformations-descent}. Next, we would like to show that 
\begin{equation}\label{equation-T-composition}
T_{(f_1,\tau_1)^*, (f_2,\tau_2)^*} \circ T_{(f_0,\tau_0)^*, (f_1,\tau_1)^*} = T_{(f_0,\tau_0)^*, (f_2,\tau_2)^*}. 
\end{equation}
Again, by considering $\tilde{S}=\tilde{S}_2\wotimes_S \tilde{S}_1 \wotimes_S \tilde{S}_0$, we may reduce to the case where the $\tau_i$ have the same target. 
Let $\mu_{(i,j)}: \tilde{R}^{(2)}\xr{} \tilde{S}$, for $(i,j)= (0,1), (1,2), (0,2)$,  be such that $\mu_{(i,j)}\circ \delta^1_0 = \tau_i$ and $\mu_{(i,j)}\circ \delta^1_1 = \tau_j$. Let $(M,\epsilon)$ be a $q$-HPD stratifications. Showing \eqref{equation-T-composition} for $M$ is equivalent to showing 
$$
\epsilon\wotimes_{\tilde{R}^{(2)}, \mu_{(1,2)}} \tilde{S} \circ \epsilon\wotimes_{\tilde{R}^{(2)}, \mu_{(0,1)}} \tilde{S} = \epsilon\wotimes_{\tilde{R}^{(2)}, \mu_{(0,2)}} \tilde{S}.
$$
To prove this we define  $\rho: \tilde{R}^{(3)}\xr{} \tilde{S}$ with $\rho\circ \delta^2_0=\mu_{(0,1)}$,  $\rho\circ \delta^2_1=\mu_{(0,2)}$, and $\rho\circ \delta^2_2 = \mu_{(1,2)}$. Applying $\wotimes_{\tilde{R}^{(3)}, \rho} \tilde{S}$ to the cocycle condition implies the claim.

We will also need the compatibility of the isomorphisms \eqref{equation-natural-transformations-descent} with base change. Let $u:S\xr{} P$ be a morphism in $(R/D)_{\qc}$, we write $\tilde{P}=\tilde{S}\wotimes_S P$ and denote by $u': \tilde{S}\xr{} \tilde{P}$ and $f':P\xr{} \tilde{P}$ the base change of $u$ and $f$, respectively. Let $u^*$ be the base change functor $\wotimes_S P$. We can identify $u^*\circ (f,\tau)^*$ with $(f',u'\circ \tau)^*$. The equality 
\begin{equation}\label{equation-T-base-change}
u^*T_{(f_0,\tau_0)^*, (f_1,\tau_1)^*} = T_{u^*\circ (f_0,\tau_0)^*, u^*\circ (f_1,\tau_1)^*}
\end{equation}
follows immediately from the construction of $T$.   
  
Finally, we define a quasi-inverse to the functor \eqref{equation-crystals-to-connections}. For each object $S$ in $(R/D)_{\qc}$ make a choice $(f_S,\tau_S)$. For a $q$-HPD stratification $(M,\epsilon)$ we set $M_S=(f_S,\tau_S)^*(M)$. For each $u:P\xr{} S$ we define $M_u:u^*M_S\xr{} M_P$ by $M_u:=T_{u^*\circ (f_S,\tau_S)^*,(f_P,\tau_P)^*}$. By using \eqref{equation-T-composition} and \eqref{equation-T-base-change} this defines a functor to $q$-crystals. For a different choice of $((f_S,\tau_S))_S$ we can use $T$ to construct a natural isomorphism between the functors. 

For a $q$-crystal $M$, the isomorphism 
$$
\tilde{S}\wotimes_{f,S} M_S \xr{M_f} M_{\tilde{S}} \xr{M_{\tau}^{-1}} \tilde{S}\wotimes_{\tau,\tilde{R}}M_{\tilde{R}}
$$   
is compatible with the trivial descent datum on $\tilde{S}\wotimes_S M_S$ and the descent datum induced by the associated $q$-HPD stratification $(M_{\tilde{R}}, M_{\delta^1_1}^{-1} \circ M_{\delta^1_0}))$ via $\tau$. It induces a natural isomorphism. 

On the other hand, starting with $q$-HPD stratifications, we may simplify the situation by choosing $\tilde{R}=\tilde{S}$, $f_{\tilde{R}}={\rm id}$, and $\tau={\rm id}_{\tilde{R}}$. For a $q$-HPD stratification $(M,\epsilon)$ with associated $q$-crystal $(S\mapsto M_S, u\mapsto M_u)$, we get $M=M_{\tilde{R}}$ and $\epsilon=M_{\delta^1_1}^{-1} \circ M_{\delta^1_0}$. 

This shows the equivalence of categories.
\end{proof} 
\end{thm}

\begin{proposition}\label{proposition-equivalence-flat}
The equivalence of categories of Theorem \ref{theorem-equivalence} induces an equivalence between the $(p,\dd)$-completely flat objects.
\begin{proof}
We only need to show that $M$ is a $(p,\dd)$-completely flat $q$-crystal if and only if $M_{\tilde{R}}$ is a  $(p,\dd)$-completely flat $\tilde{R}$-module. Then use Proposition \ref{proposition-retractions} together with Lemma \ref{lemma-base-change-flatness}, Lemma \ref{lemma-tensor-product-flatness}, and Lemma \ref{lemma-descent-flatness}.
\end{proof}
\end{proposition}

\section{$q$-connections}

In this section, we let $\tilde{R}$ be as in Proposition \ref{proposition-retractions}. We will use the notation from \textsection\ref{subsection-stratifiactions}. This section will borrow from \cite[\textsection4]{BO}.  

\begin{definition}
For $(p,\dd)$-derived complete $\tilde{R}$-modules $M$ and $N$, we define the {\em $q$-HPD differential operators} from $M$ to $N$ by
$$
q{\rm -HPDiff}_{\tilde{R}}(M,N)=\Hom_{\tilde{R}}((M\wotimes_{\tilde{R},\delta^1_0} \tilde{R}^{(2)})_{\delta^1_1}, N),
$$ 
where $(.)_{\delta^1_1}$ means that we consider it as an $\tilde{R}$-module via $\delta^1_1$.
\end{definition}

We would like to make $(p,\dd)$-derived complete $\tilde{R}$-modules together with $q$-HPD differential operators as morphisms into a category. In order to define the composition, we set 
$$
V= \tilde{R}^{(2)}\wotimes_{\delta^1_1,\tilde{R},\delta^1_0} \tilde{R}^{(2)},
$$
which is an object in $(R/D)_{\qc}$, and note that there is a unique morphism $\theta:\tilde{R}^{(3)}\xr{} V$ in $(R/D)_{\qc}$ such that $\theta\circ \delta^2_0={\rm id}_{\tilde{R}^{(2)}}\wotimes 1$ and $\theta\circ \delta^2_2=1 \wotimes {\rm id}_{\tilde{R}^{(2)}}$. 

Now, let $f\in q{\rm -HPDiff}_{\tilde{R}}(M,N)$ and $g\in q{\rm -HPDiff}_{\tilde{R}}(N,P)$, we can form the composition 
\begin{equation}\label{equation-pre-composition}
M\wotimes_{\tR,\delta^1_0} \tR^{(2)}\xr{{\rm id}_M\wotimes \theta\circ \delta^2_1} M\wotimes_{\tilde{R},\delta^1_0} \tR^{(2)}\wotimes_{\delta^1_1,\tilde{R},\delta^1_0} \tilde{R}^{(2)} \xr{f\wotimes {\rm id}_{\tR^{(2)}}} N\wotimes_{\tilde{R},\delta^1_0} \tR^{(2)} \xr{g} P, 
\end{equation}
which is $\tR$-linear if we consider the source as an $\tR$-module via $\delta^1_1$. We define $g\circ f$ as the composition \eqref{equation-pre-composition}. 

In order to define ${\rm id}_M\in q{\rm -HPDiff}_{\tilde{R}}(M,M)$, we note that there is a unique morphism ${\rm mult}:\tR^{(2)}\xr{} \tR$ in $(R/D)_{\qc}$ such that ${\rm mult}\circ \delta^1_0 = {\rm id}_{\tilde{R}} = {\rm mult}\circ \delta^1_1$.  We take 
$$
(M\wotimes_{\tilde{R},\delta^1_0} \tilde{R}^{(2)})_{\delta^1_1} \xr{{\rm id}_M\wotimes {\rm mult}} M
$$  
as the identity. This finishes the definition of the category $q{\rm -HPDiff}_{\tilde{R}}$. 

Next, we define a functor 
\begin{equation}\label{equation-to-D-modules}
q{\rm -HPDiff}_{\tilde{R}} \xr{} {\rm Mod}^{\wedge}_D,
\end{equation}
which takes a $(p,\dd)$-derived complete $\tilde{R}$-module to a $(p,\dd)$-derived complete $D$-module by restriction via $D\xr{} \tilde{R}$.  For $f\in q{\rm -HPDiff}_{\tilde{R}}(M,N)$ we simply define 
$$
M\xr{m\mapsto m\wotimes 1} M\wotimes_{\tilde{R},\delta^1_0} \tilde{R}^{(2)}\xr{f} N
$$
as the corresponding morphism. 

The functor defined in \eqref{equation-to-D-modules} is not faithful (see \cite[\textsection4.4]{BO}). However, the map 
\begin{equation}\label{equation-HPD-Hom-to-D-Hom}
q{\rm -HPDiff}_{\tilde{R}}(M,N) \xr{} \Hom_{D}(M,N)
\end{equation}  
is injective if $N$ is $\dd$-torsion free. 

\subsection{}
Let $(M,\epsilon)$ be a $q$-HPD stratification on $(\tilde{R}, I\tilde{R})$. To simplify the notation, we set 
$
\qdiff(N):=\qdiff(N,N).
$

By using 
\begin{align*}
M\wotimes_{\tR,\delta^1_1} \tR^{(2)} & \xr{} \Hom_{\tR}(\qdiff(\tR), M)\\
m\otimes r &\mapsto [\xi \mapsto \xi(r)\cdot m],  
\end{align*}
we get a map
\begin{equation}\label{equation-epsilon-to-nabla}
\Hom_{\tR^{(2)}}(M\wotimes_{\tR,\delta^1_0} \tR^{(2)}, M\wotimes_{\tR,\delta^1_1} \tR^{(2)}) \xr{} \Hom_{\tR}(\qdiff(\tR), \qdiff(M)).
\end{equation}
We denote by $\nabla$ the image of $\epsilon$. Explicitly, $\nabla(\xi)(m\otimes t) = ({\rm id}_M\wotimes \xi)(\epsilon(m\otimes t))$. 

\begin{lemma}\label{lemma-cocycle-condition-composition}
The map $\nabla$ respects compositions. That is, for any $\zeta,\xi \in \qdiff(\tR)$ we have $\nabla(\zeta\circ \xi)=\nabla(\zeta)\circ \nabla(\xi)$.
\begin{proof}
Recall that $V= \tilde{R}^{(2)}\wotimes_{\delta^1_1,\tilde{R},\delta^1_0} \tilde{R}^{(2)}$ and we have $\theta:\tR^{(3)}\xr{} V$. As in \eqref{equation-pre-composition}, we consider 
\begin{equation}
M\wotimes_{\tilde{R},\delta^1_0} \tR^{(2)}\wotimes_{\delta^1_1,\tilde{R},\delta^1_0} \tilde{R}^{(2)} \xr{\nabla(\xi)\wotimes {\rm id}_{\tR^{(2)}}} M\wotimes_{\tilde{R},\delta^1_0} \tR^{(2)} \xr{\nabla(\zeta)} M. 
\end{equation}
This map can be rewritten as follows
\begin{multline}\label{multline-first-step-compatibility}
M\wotimes_{\tR, \theta \circ  \delta^2_0\circ \delta^1_0} V \xr{(\theta\circ \delta^2_0)^*(\epsilon)} M\wotimes_{\tR, \theta \circ  \delta^2_2\circ \delta^1_0} V \xr{(\theta\circ \delta^2_2)^*(\epsilon)} M\wotimes_{\tR, \theta \circ  \delta^2_2\circ \delta^1_2} V \\ \xr{{\rm id}_M\wotimes \zeta\circ (\xi\wotimes {\rm id}_{\tR^{(2)}})} M,  
\end{multline}
where $f^*(\epsilon):=\epsilon\wotimes_{\tR^{(2)}, f} V$, and the last arrow is induced by the last two arrows in \eqref{equation-pre-composition} for $f=\xi$ and $g=\zeta$. 

We use the cocycle condition to identify the composition of the first two arrows in \eqref{multline-first-step-compatibility} with $(\theta\circ\delta^2_1)^*(\epsilon)$. Then the claim follows from the commutativity of the following diagram:
$$
\xymatrix@C+=3cm{
M\wotimes_{\tR,\delta^1_0} \tR^{(2)} \ar[r]^{\epsilon} \ar[dd]_{{\rm id}_M\wotimes \theta\circ \delta^2_1} &  M\wotimes_{\tR,\delta^1_1} \tR^{(2)} \ar[dd]_{{\rm id}_M\wotimes \theta\circ \delta^2_1} \ar[rd]^-{{\rm id}_M\wotimes (\zeta\circ \xi)}
\\
& & M
\\
M\wotimes_{\tR, \theta \circ  \delta^2_0\circ \delta^1_0} V \ar[r]^{(\theta\circ\delta^2_1)^*(\epsilon)} &  M\wotimes_{\tR, \theta \circ  \delta^2_2\circ \delta^1_2} V \ar[ru]_-{\quad\quad {\rm id}_M\wotimes \zeta\circ (\xi\wotimes {\rm id}_{\tR^{(2)}})} 
}
$$
\end{proof}
\end{lemma}

\newcommand{\qconn}{q{\rm -Conn}}
\begin{definition}\label{definition-q-connections}
A {\em $q$-connection} on $\tilde{R}$ is a $(p,\dd)$-derived complete $D$-module $M$ together with morphism of $D$-algebras $\nabla\in \Hom_D(\qdiff(\tR), \Hom_D(M,M))$.  

Morphisms of $q$-connections are morphisms of $D$-modules that are compatible with $\nabla$ in the obvious way.
This defines the {\em category of $q$-connections} $\qconn_{(\tilde{R},I)}$.
 
Naturally, every $q$-connection is an $\tilde{R}$-module and every morphism of $q$-connections is a morphism of $\tilde{R}$-modules. A $q$-connection is called $(p,\dd)$-completely flat if $M$ is a $(p,\dd)$-completely flat $\tilde{R}$-module. 
\end{definition}

Lemma \ref{lemma-cocycle-condition-composition} shows that $(M,\epsilon)\mapsto (M,\nabla)$ from \eqref{equation-epsilon-to-nabla} (and using \eqref{equation-to-D-modules}) yields a functor 
\begin{equation}\label{equation-stratification-to-connections}
\text{($q$-HPD stratifications on $(\tilde{R}, I\tilde{R})$)} \xr{} \qconn_{(\tilde{R},I)}.
\end{equation}
Our next goal is to show an equivalence between the full subcategories of $(p,
\dd)$-completely flat $q$-HPD stratifications and quasi-nilpotent $(p,
\dd)$-completely flat $q$-connections induced by this functor. We will start by finding a nice basis for $\tilde{R}^{(2)}$ lifting the standard basis (given by the coordinates) modulo $q-1$. 

Recall that we have $\tilde{R}\lf \epsilon_1,\dots,\epsilon_n\rf \xr{} \tilde{R}^{(2)}$. For $I\in \Z_{\geq 0}^n$ we set $\epsilon^I:=\prod_{k=1}^n\epsilon_k^{I_k}$, and $[I!]_p$ is defined as the largest $p$-power dividing $\prod_{k=1}^n I_k!$. 

As in the proof of Proposition \ref{proposition-retractions}, we can find a sequence $(\Gamma_I)_{I\in \Z_{\geq 0}^n}$ of elements in $\tilde{R}^{(2)}$ such that $[I!]_p \cdot \Gamma_I \equiv \epsilon^I \mod (q-1)$.
%\begin{definition}
%A sequence $(\Gamma_I)_{I\in \Z_{\geq 0}^n}$ of elements in $\tilde{R}^{(2)}$ is called an {\em $(q-1)$-adically converging lift} if the following conditions are satisfied:
%\begin{enumerate}
%\item For all $I\in \Z_{\geq 0}^n$ there is $n_I\in D$ such that $n_I\equiv [I!]_p \mod q-1$, $n_I\cdot \Gamma_I \in \tilde{R}\lf \epsilon_1,\dots,\epsilon_n\rf$, and  $n_I\cdot \Gamma_I \equiv \epsilon^I \mod (q-1)$.
%\item We have $\lim_I n_I\cdot \Gamma_I=0$ in the $(\epsilon_1,\dots,\epsilon_n,q-1)$-adic sense. In other words, for each $k$, there are only finitely many indices $I$ such that $$n_I\cdot \Gamma_I \not\in (\epsilon_1,\dots,\epsilon_n,q-1)^k.$$ 
%\end{enumerate}
%\end{definition}
And we obtain an isomorphism 
$$
\left( \bigoplus_{I\in \Z_{\geq 0}^n} \delta^1_1(\tilde{R}) \cdot \Gamma_I  \right)^{\wedge} \xr{\cong} \tilde{R}^{(2)},
$$
where the source is the derived $(p,\dd)$-completion of the direct sum $\bigoplus_{I\in \Z_{\geq 0}^n}$, but turns out to be the classical $(p,\dd)$-completion, and is automatically the direct sum in the category of $(p,\dd)$-derived complete modules. In particular, if $(\xi_I)_{I\in \Z_{\geq 0}^n}$ denotes the dual basis then 
$$
\qdiff(\tilde{R}) = \prod_{I\in \Z_{\geq 0}^n} \tilde{R}\cdot \xi_I.
$$

%\begin{lemma}\label{lemma-existence-basis}
%A $(q-1)$-adically converging lift $(\Gamma_I)_{I\in \Z_{\geq 0}^n}$ exists. 
%\begin{proof}
%Let $K= (\epsilon_1,\dots,\epsilon_n,q-1)$, as ideal in the 
%$\delta$-$\tilde{R}$-algebra $\tilde{R}\lf \epsilon_1,\dots,\epsilon_n\rf$ used in the construction of $\tilde{R}^{(2)}$. We have $\delta(K)\subset K$, hence $\delta(K^m)\subset K^m$. Then we can simply define $(\Gamma_I)_{I\in \Z_{\geq 0}^n}$ as in the proof of Proposition \ref{proposition-retractions}.
%\end{proof}
%\end{lemma}

%We will fix a $(q-1)$-adically converging lift $(\Gamma_I)_{I\in \Z_{\geq 0}^n}$ %with dual basis $(\xi_I)_{I\in \Z_{\geq 0}^n}$ in the following.

%\begin{lemma}
%For all $K$ and $t\in \tilde{R}^{(2)}$, we have $\lim_I \xi_K(\Gamma_I\cdot t)=0$ in the $(p,\dd)$-adic sense. 
%\begin{proof}
%We can reduce to the case $t=\Gamma_J$. In this case, we prove the stronger convergence in the $(q-1)$-adic sense. By using that $\tilde{R}/(q-1)\tilde{R}$ is  $p$-torsionfree we may reduce to proving  $\lim_I \xi_K(n_I\Gamma_I\cdot n_J\Gamma_J)=0$. It is not hard to finish the proof from here, but it is a bit ugly.
%\end{proof}
%\end{lemma}

\begin{definition}\label{definition-quasi-nilpotence}
A $(p,\dd)$-completely flat $q$-connection $(M,\nabla)$ is called {\em quasi-nilpotent} if for each $m\in M$, we have $\lim_I \nabla(\xi_I)(m)=0$ in the $(p,\dd)$-adic sense. 
\end{definition}

\begin{proposition}\label{proposition-equivalence-stratification-connections}
The functor \eqref{equation-stratification-to-connections} 
%from $q$-HPD stratifications to $q$-connections 
induces an equivalence between $(p,\dd)$-completely flat $q$-HPD stratifications and $(p,\dd)$-completely flat quasi-nilpotent $q$-connections.
\begin{proof}
Let $M$ be a $(p,\dd)$-derived complete and $(p,\dd)$-completely flat $\tilde{R}$-module. Then 
$$
\left( \bigoplus_{I\in \Z_{\geq 0}^n} M \otimes \Gamma_I  \right)^{\wedge} \xr{\cong}  M\wotimes_{\tilde{R}, \delta^1_1} \tilde{R}^{(2)}
$$
is an isomorphism and the left hand side is the same as the classical completion. Therefore  
$$
 M\wotimes_{\tilde{R}, \delta^1_1} \tilde{R}^{(2)} \xr{(\xi_I)_I} \prod_{I\in \Z_{\geq 0}^n} M
$$
is injective and the image equals
$
\{(m_I)_I \mid \lim_I m_I=0\},
$
where the limit is in the $(p,\dd)$-adic sense.

Now let $(M,\epsilon)$ be a $q$-HPD stratification with associated $q$-connection $(M,\nabla)$. By definition, $\nabla(\xi_I)(m)=\xi_I(\epsilon(m))$, which implies $\lim_I \nabla(\xi_I)(m)=0$, hence $(M,\nabla)$ is quasi-nilpotent. 

Given a quasi-nilpotent $(M,\nabla)$, we can construct $\epsilon$ as follows. We want
$$
\epsilon(m\otimes 1) = \sum_I \nabla(\xi_I)(m)\otimes \Gamma_I,
$$  
which will give a well-defined morphism after showing that 
$$
\sum_I \nabla(\xi_I)(a\cdot m)\otimes \Gamma_I = \sum_I \nabla(\xi_I)(m)\otimes \Gamma_I\cdot \delta^1_0(a)
$$
for all $a\in \tilde{R}$. This is equivalent to 
$$
\xi_K(a\cdot m) = \sum_I \xi_K(\Gamma_I\cdot \delta^1_0(a))\cdot \nabla(\xi_I)(m)
$$
for all $K$, and follows from 
$$
\xi_K\circ a = \sum_I \xi_K(\Gamma_I\cdot \delta^1_1(a)) \cdot \xi_I 
$$
in $\qdiff(\tilde{R})$. But this just means 
$$
(\xi_K\circ a)(\Gamma_I) = \xi_K(\Gamma_I\cdot \delta^1_1(a))
$$
for all $I$, and holds by definition of the composition in $\qdiff(\tilde{R})$. 

Now that we have constructed $\epsilon$, we need to show the cocycle condition. One easily computes 
\begin{align*}
(\delta^{2*}_2(\epsilon) \circ \delta^{2*}_0(\epsilon))(m\otimes 1) &= \sum_{I,J} \nabla(\xi_J)(\nabla(\xi_I)(m))\otimes \delta^2_2(\Gamma_J)\cdot \delta^2_0(\Gamma_I)\\
\delta^{2*}_1(\epsilon)(m\otimes 1) &= \sum_K \nabla(\xi_K)(m)\otimes \delta^2_1(\Gamma_K).   
\end{align*}
We have 
$$
\left( \bigoplus_{I,J} \delta^2_1(\tilde{R}) \cdot \delta_0^2(\Gamma_I)\delta^2_2(\Gamma_J)  \right)^{\wedge} \xr{\cong} \tilde{R}^{(3)},
$$
and the "dual" basis is given by 
$$
\tilde{R}^{(3)} \xr{\theta} V \xr{\xi_I\wotimes {\rm id}_{\tilde{R}^{(2)}}} \tilde{R}^{(2)} \xr{\xi_J} \tilde{R}.  
$$ 
Let us write 
$$
\delta^2_1(\Gamma_K) = \sum_{I,J} \delta^2_1(t_{I,J}(K)) \cdot \delta_0^2(\Gamma_I)\delta^2_2(\Gamma_J),
$$
which automaticall implies $\lim_{I,J} t_{I,J}(K) = 0$. Then we get 
$$
\xi_J \circ \xi_I = \sum_K t_{I,J}(K) \cdot \xi_K
$$
in $\qdiff(\tilde{R})$ by definition of the composition. This proves
$$
(\delta^{2*}_2(\epsilon) \circ \delta^{2*}_0(\epsilon))(m\otimes 1) = \delta^{2*}_1(\epsilon)(m\otimes 1),
$$ 
and the cocycle condition.

The functors are inverse to each other on the nose. 
\end{proof} 
\end{proposition}

\subsection{}
Our next goal is to understand how the categories of $q$-connections for two different liftings $\tilde{R}_1$ and $\tilde{R}_2$ are related. Theorem \ref{theorem-equivalence} and Proposition \ref{proposition-equivalence-stratification-connections} tell us that the full subcategories of $(p,\dd)$-completely flat and quasi-nilpotent $q$-connections are equivalent. 

Recall that we work in $(R/D)_{\qc}$, with base the $q$-PD pair $(D,I)$. With lifting we mean $(\tilde{R}_i,I\tilde{R}_i)$ are objects in $(R/D)_{\qc}$ such that $\tilde{R}_i/I\tilde{R}_i=R$. In particular, $\tilde{R}_i$ is a $D$-$\delta$-algebra. 

By using the arguments of Proposition \ref{proposition-retractions}, one can find an isomorphism of $D$-algebras $\tau':\tilde{R}_1\xr{} \tilde{R}_2$ inducing the identity modulo $I$. However, this isomorphism is not compatible with the $\delta$-structures, hence not a morphism in $(R/D)_{\qc}$. In general, finding an isomorphism in $(R/D)_{\qc}$ is not possible.  

Let us consider the classical crystalline situation $q-1=0$ for a moment. This simplifies the situation, because we can drop $\delta$ and use PD-ideals only.  In this case, $\tilde{R}_i$ is $p$-torsion free and $I\tilde{R}_i$ is a PD-ideal. Moreover, $\qdiff(\tilde{R}_i)$ becomes the non-commutative $\tilde{R}_i$-algebra formally generated by $\partial_{x_1},\dots,\partial_{x_n}$ for some coordinates $x_1,\dots,x_n$ (whose existence we assume). Explicitly, 
$$
\qdiff(\tilde{R}_i) = \prod_{I\in \Z^n_{\geq 0}} \tilde{R}_i\cdot \partial_{x_1}^{I_1}\cdots \partial_{x_n}^{I_n},
$$
(to see the independence of the choice of coordinates, we have to use that $\tilde{R}_i$ is classically $p$-complete.) In this case, we can simply use $\tau'$ to construct an isomorphism of $D$-algebras
$$
\tau':\qdiff(\tilde{R}_1) \xr{\cong}  \qdiff(\tilde{R}_2),
$$
which will induce an equivalence of categories $(\tau')^*$ between the $q$-connections on $\tilde{R}_1$ and $\tilde{R}_2$, respectively. For another choice of an isomorphism $\tau''$, there is a natural isomorphism $(\tau')^*\xr{} (\tau'')^*$

Let us go back to the general case (where maybe $q-1\neq 0$). Then $\qdiff(\tilde{R}_i)$ depends on the $\delta$-structure of $\tilde{R}_i$, and $\tau'$ cannot be used.  However, it is possible to use the construction from Proposition \ref{proposition-retractions} to define an isomorphism of $D$-algebras
$$
s:\qdiff(\tilde{R}_1) \xr{\cong}  \qdiff(\tilde{R}_2),
$$ 
inducing an equivalence of categories between the categies of $q$-connections. Again, $s$ is not unique, but two choices are naturally isomorphic.   When restricted to $(p,\dd)$-completely flat and quasi-nilpotent $q$-connections this equivalence is compatible with the equivalence from Theorem \ref{theorem-equivalence}. 

\subsection{} \label{subsection-compare-two-liftings}
Let us recall the constructions from Proposition \ref{proposition-retractions} for   $\tilde{R}=\tilde{R}_1$, $S=\tilde{R}_2$. We obtain 
$$
\xymatrix{
& \tilde{S} & 
\\
\tilde{R}_1 \ar@{.>}[rr]^{\tau'} \ar[ur]^{\tau} & & \tilde{R}_2 \ar[ul]_{f} 
}
$$
where the dotted arrow is an isomorphism, but is not compatible with the $\delta$-structure. Moreover, $\tilde{S}$ is the $q$-PD-envelope of $\tilde{R}_2\lf \epsilon_1,\dots,\epsilon_n \rf$ for the regular sequence $\epsilon_1,\dots,\epsilon_n$, and $\tau$ and $f$ factor over $\tilde{R}_2\lf \epsilon_1,\dots,\epsilon_n \rf$. Explicitly, we have $\tau(x_i)=\tau'(x_i)+\epsilon_i$, for coordinates $x_1,\dots,x_n$, and $f(a)=a$. However, $\tilde{S}$ does not depend on the coordinates $(x_i)$. Even better, it does not depend on $\tau'$. It is the coproduct of $\tilde{R}_1$ and $\tilde{R}_1$ in the category $(R/D)_{\qc}$. We will use the notation $\tilde{R}_1 \otimes_{\qc} \tilde{R}_2:=\tilde{S}$, $\delta_0:=\tau$, and $\delta_1:=f$.

Note that $\tilde{R} \otimes_{\qc} \tilde{R}=\tilde{R}^{(2)}$ via the morphisms $\delta^1_0$ and $\delta^1_1$. In the following, we will modify the constructions which were used to define the algebra structure on $\qdiff(\tilde{R})$ in order to obtain 
\begin{multline}\label{multline-composition-U}
\Hom_{\tilde{R}_3}(\tilde{R}_2 \otimes_{\qc} \tilde{R}_3,\tilde{R}_3) \times \Hom_{\tilde{R}_2}(\tilde{R}_1 \otimes_{\qc}\tilde{R}_2,\tilde{R}_2) \xr{}\\ \Hom_{\tilde{R}_3}(\tilde{R}_1 \otimes_{\qc}\tilde{R}_3,\tilde{R}_3),
\end{multline}
which we will write as composition $(t,s)\mapsto t\circ s$.

By using the universal property we get 
$$
\rho: \tilde{R}_1 \otimes_{\qc} \tilde{R}_3 \xr{} (\tilde{R}_1\otimes_{\qc}\tilde{R}_2) \wotimes_{\delta_1, \tilde{R}_2, \delta_0} (\tilde{R}_2\otimes_{\qc}\tilde{R}_3)
$$
such that $\rho\circ \delta_0 = \delta_0\wotimes 1$ and $\rho\circ \delta_1=1\wotimes \delta_1$. The composition \eqref{multline-composition-U} is defined by 
\begin{multline*}
\tilde{R}_1 \otimes_{\qc}\tilde{R}_3 \xr{\rho} (\tilde{R}_1\otimes_{\qc}\tilde{R}_2) \wotimes_{\delta_1, \tilde{R}_2, \delta_0} (\tilde{R}_2\otimes_{\qc}\tilde{R}_3) \xr{s\wotimes {\rm id}} \\ \tilde{R}_2\otimes_{\qc}\tilde{R}_3 \xr{t} \tilde{R}_3.
\end{multline*}
It is associative and compatible with the composition on $\qdiff(\tilde{R})$ introduced in \eqref{equation-pre-composition}. In particular, $\Hom_{\tilde{R}_2}(\tilde{R}_1\otimes_{\qc}\tilde{R}_2,\tilde{R}_2)$ is a right $\qdiff(\tilde{R}_1)$-module and a left $\qdiff(\tilde{R}_2)$-module.

We define 
\begin{multline*}
{\rm SHom}_{\tilde{R}_2}(\tilde{R}_1\otimes_{\qc}\tilde{R}_2,\tilde{R}_2) =\{s\in  {\rm Hom}_{\tilde{R}_2}(\tilde{R}_1\otimes_{\qc}\tilde{R}_2 ,\tilde{R}_2) \mid \\ s(1)=1, (s\circ \delta_0)\otimes_D D/I = {\rm id}_R\}.
\end{multline*}
Note that $s(1)=1$ implies $s\circ \delta_1={\rm id}_{\tilde{R}_2}$, in other words, $s$ is a section. 

We know that ${\rm SHom}_{\tilde{R}_2}(\tilde{R}_1\otimes_{\qc}\tilde{R}_2,\tilde{R}_2)$ is not empty, because the section constructed in the proof of Proposititon \ref{proposition-retractions} satisfies the requirements (indeed, it satisfies $s\circ \delta_0 \otimes_{D} D/(q-1) = \tau'\otimes_{D} D/(q-1)$, and $\tau'\otimes_{D} D/I = {\rm id}_{R}$ holds by definition).  

The composition \eqref{multline-composition-U} respects ${\rm SHom}$, because $(t\circ s)\circ \delta_0 = (t\circ \delta_0) \circ (s\circ \delta_0)$.  
Moreover,  ${\rm SHom}_{\tilde{R}}(\tilde{R}\otimes_{\qc}\tilde{R} ,\tilde{R})$ is a multplicative subgroup of $\qdiff(\tilde{R})$, because $I^p\tilde{R}\subset (p,\dd)\tilde{R}$ and the injectitivity of the map \eqref{equation-HPD-Hom-to-D-Hom} 
$$\qdiff(\tilde{R}) \xr{} \Hom_D(\tilde{R},\tilde{R}), \quad \xi\mapsto \xi\circ \delta^1_0.$$ 

\begin{lemma}\label{lemma-existence-inverse}
For $s\in {\rm SHom}_{\tilde{R}_2}(\tilde{R}_1\otimes_{\qc}\tilde{R}_2,\tilde{R}_2)$ there exists a unique $t\in {\rm SHom}_{\tilde{R}_1}(\tilde{R}_2\otimes_{\qc}\tilde{R}_1,\tilde{R}_1)$ such that $s\circ t=1$ and $t\circ s=1$.
\begin{proof}
We know that ${\rm SHom}_{\tilde{R}_1}(\tilde{R}_2\otimes_{\qc}\tilde{R}_1,\tilde{R}_1)$ is non empty. So we can find some $t'$. Then $t'\circ s\in {\rm SHom}_{\tilde{R}_1}(\tilde{R}_1\otimes_{\qc}\tilde{R}_1 ,\tilde{R}_1)$, and we set
$t= (t'\circ s)^{-1}\circ t'$. Now that we have $t\circ s=1$, we play the same game with $t$ to get $u$ with $u\circ t=1$, hence $u=u\circ t\circ s=s$.  
\end{proof}
\end{lemma}

\begin{proposition}\label{proposition-morphism-qdiff}
Every $s\in {\rm SHom}_{\tilde{R}_2}(\tilde{R}_1\otimes_{\qc}\tilde{R}_2,\tilde{R}_2)$ induces an isomorphism 
$$
\psi_s: \qdiff(\tilde{R}_1) \xr{} \qdiff(\tilde{R}_2). 
$$
For a second choice $s'\in {\rm SHom}_{\tilde{R}_2}(\tilde{R}_1\otimes_{\qc}\tilde{R}_2,\tilde{R}_2)$ there is a $\xi\in {\rm SHom}_{\tilde{R}_2}(\tilde{R}_2\otimes_{\qc}\tilde{R}_2 ,\tilde{R}_2)$ such that $\psi_{s'}(\zeta)=\xi \cdot \psi_{s}(\zeta)\cdot \xi^{-1}$ for all $\zeta$. In other words, $\psi_s$ is up to conjugation by elements in ${\rm SHom}_{\tilde{R}_2}(\tilde{R}_2\otimes_{\qc}\tilde{R}_2 ,\tilde{R}_2)$ independent of the choice of $s$.
\begin{proof}
We can use Lemma \ref{lemma-existence-inverse} to find $t$, and define $\psi_{s}(\zeta)=s\circ \zeta\circ t$ \eqref{multline-composition-U}. 

Making another choice $s'$, we define $\xi:=s'\circ t$. Then $s'=\xi\circ s$ and $\psi_{s'}(\zeta)=\xi \cdot \psi_{s}(\zeta)\cdot \xi^{-1}$ for all $\zeta$.
\end{proof}
\end{proposition}

\begin{corollary}\label{corollary-equivalance-q-connections}
Let $(\tilde{R}_i, I\tilde{R}_i)$ in $(R/D)_{\qc}$, for $i=0,1,2,3$, be  such that $\tilde{R}_i/I\tilde{R}_i=R$, and suppose that $(p, \dd)$-\'etale maps $h_i:D[x_1,\dots,x_n]^{\wedge}\xr{} {R}_i$ exist. 

For every $s_{2,1}\in {\rm SHom}_{\tilde{R}_2}(\tilde{R}_1\otimes_{\qc}\tilde{R}_2,\tilde{R}_2)$ we have an equivalence of categories 
\begin{align*}
F_{s_{2,1}}: \qconn_{(\tilde{R}_2,I\tilde{R}_2)} &\xr{}  \qconn_{(\tilde{R}_1,I\tilde{R}_1)}\\
(M,\nabla) &\mapsto (M,\nabla\circ \psi_s).
\end{align*}
Between two equivalences of this form there is a natural isomorphism:
$$
u_{\xi}: F_s \xr{} F_{s'}, \quad u_{\xi, (M,\nabla)}(m):=\nabla(\xi)(m),
$$ 
where $\xi$ is defined by $s'=\xi\circ s$.

For $s_{3,1}\in {\rm SHom}_{\tilde{R}_3}(\tilde{R}_1\otimes_{\qc}\tilde{R}_3,\tilde{R}_3)$ and $s_{3,2}\in {\rm SHom}_{\tilde{R}_3}(\tilde{R}_2\otimes_{\qc}\tilde{R}_3,\tilde{R}_3)$, we have $F_{s_{2,1}} \circ F_{s_{3,2}} = F_{s_{3,2}\circ s_{2,1}}$. In particular, we have a natural isomorphism
$$
t_{3,2,1}: F_{s_{2,1}} \circ F_{s_{3,2}} = F_{s_{3,2}\circ s_{2,1}} \xr{u_{\zeta}} F_{s_{3,1}},
$$ 
where  $\zeta$ is such that $s_{3,1}=\zeta \circ s_{3,2}\circ s_{2,1}$.

Moreover, for all $(M,\nabla)$ the following diagram is commutative:
$$
\xymatrix{
F_{s_{1,0}} (F_{s_{2,1}} \circ F_{s_{3,2}}(M,\nabla)) \ar[rr]^{F_{s_{1,0}}(t_{3,2,1})} \ar[d]_{t_{2,1,0}(F_{s_{3,2}}(M,\nabla))} & &
F_{s_{1,0}} ( F_{s_{3,1}}(M,\nabla) ) \ar[d]^{t_{3,1,0}} \\
F_{s_{2,0}}(F_{s_{3,2}}(M,\nabla)) \ar[rr]^{t_{3,2,0}} & &
F_{s_{3,0}}(M,\nabla).
}
$$
\begin{proof}
This follows immediately from Proposition \ref{proposition-morphism-qdiff}. 
\end{proof}
\end{corollary}

\begin{proposition}
Assumptions as in Corollary \ref{corollary-equivalance-q-connections}. The following diagram of functors commutes up to natural transformation
$$
\xymatrix{
& q{\rm -Cris}(R/D) \ar[dr]^{\text{\eqref{equation-crystals-to-connections}}} \ar[dl]_{\text{\eqref{equation-crystals-to-connections}}} & \\
\text{($q$-HPD str.~on $(\tilde{R}_2, I\tilde{R}_2)$)} \ar[d]_{\text{\eqref{equation-stratification-to-connections}}} & & \text{($q$-HPD str.~on $(\tilde{R}, I\tilde{R}_1)$)} \ar[d]^{\text{\eqref{equation-stratification-to-connections}}} \\
\qconn_{(\tilde{R}_2,I\tilde{R}_2)} \ar[rr]^{\text{Corollary \ref{corollary-equivalance-q-connections}}} & & \qconn_{(\tilde{R}_1,I\tilde{R}_1)}
}
$$ 
\end{proposition}

\subsection{} 
The sections in ${\rm SHom}_{\tilde{R}_2}(\tilde{R}_1\otimes_{\qc}\tilde{R}_2,\tilde{R}_2)$ used to define the functor between the categories of $q$-connections are $p$-adic in nature. For the next section, where we will prove a global version of Corollary \ref{corollary-equivalance-q-connections}, we need to isolate those sections that come from global ones. This will not work for $p=2$.

Recall the setup from Subsection \ref{subsection-compare-two-liftings}. The $D$-algebra  $\tilde{R}_2\lf \epsilon_1,\dots,\epsilon_n \rf$ together with the ideal $I+(\epsilon_1,\dots,\epsilon_n)$ and the morphisms $f:\tilde{R}_2\xr{} \tilde{R}_2\lf \epsilon_1,\dots,\epsilon_n \rf$ and $\tau: \tilde{R}_1\xr{} \tilde{R}_2\lf \epsilon_1,\dots,\epsilon_n \rf$ does not depend on the choice of the coordinates or $\tau'$ or the $\delta$-structures. More conceptually, it is the coproduct of $\tilde{R}_1$ and $\tilde{R}_2$ in a suitable category of infinitisimal thickenings of $R$. We will write 
$$
\tilde{R}_1\otimes_{{\rm inf}} \tilde{R}_2:= \tilde{R}_2\lf \epsilon_1,\dots,\epsilon_n \rf, \quad J:= I\cdot \tilde{R}_2\lf \epsilon_1,\dots,\epsilon_n \rf + (\epsilon_1, \dots ,\epsilon_n) ,
$$    
and $\delta_1:=f$ and $\delta_0:=\tau$. Moreover, if $K\subset \tilde{R}_2$  is an ideal, then we will use the notation
\begin{align*}
{\rm Diff}^{\wedge}_{K}(\tilde{R}_1,\tilde{R}_2):= \{s\in \Hom_{\tilde{R}_2}(\tilde{R}_1\otimes_{{\rm inf}} \tilde{R}_2, \tilde{R}_2)\mid \lim_k s(J^k)=0 \quad  \text{$K$-adically} \}.
\end{align*}
Actually, we will only consider $K=(p,\dd)$ and $K=I$. Since $I^p\subset (p,\dd)$, we have an injective map
\begin{equation}\label{equation-global-to-local}
{\rm Diff}^{\wedge}_{I}(\tilde{R}_1,\tilde{R}_2) \xr{} {\rm Diff}^{\wedge}_{(p,\dd)}(\tilde{R}_1,\tilde{R}_2).
\end{equation}
As in \eqref{multline-composition-U}, we get composition maps 
$$
{\rm Diff}^{\wedge}_{K}(\tilde{R}_2,\tilde{R}_3) \times {\rm Diff}^{\wedge}_{K}(\tilde{R}_1,\tilde{R}_2) \xr {} {\rm Diff}^{\wedge}_{K}(\tilde{R}_1,\tilde{R}_3)
$$ 
for $K=(p,\dd)$ and $K=I$. This turns ${\rm Diff}^{\wedge}_{K}(\tilde{R},\tilde{R})$ into an algebra. 

We set 
\begin{equation*}
{\rm SDiff}^{\wedge}_{K}(\tilde{R}_1,\tilde{R}_2) =\{s\in {\rm Diff}^{\wedge}_{K}(\tilde{R}_1,\tilde{R}_2) \mid  s(1)=1, (s\circ \delta_0)\otimes_D D/I = {\rm id}_R\}.
\end{equation*}

\begin{proposition}
The map ${\rm Hom}_{\tilde{R}_2}(\tilde{R}_1\otimes_{\qc}\tilde{R}_2,\tilde{R}_2) \xr{} \Hom_{\tilde{R}_2}(\tilde{R}_1\otimes_{{\rm inf}} \tilde{R}_2, \tilde{R}_2)$, induced by $\tilde{R}_1\otimes_{{\rm inf}} \tilde{R}_2 \xr{} \tilde{R}_1\otimes_{\qc} \tilde{R}_2$, factors through ${\rm Diff}^{\wedge}_{(p,\dd)}(\tilde{R}_1,\tilde{R}_2)$.
\begin{proof}
This follows immediately from $J^p\subset (p,\dd)$, where $J\subset \tilde{R}_1\otimes_{\qc} \tilde{R}_2$ is the $q$-PD ideal.
\end{proof}
\end{proposition}

\begin{proposition}\label{proposition-q-1-adic-sections}
Suppose $p>2$ and $(p,\dd)$ is a regular sequence in $D$. Then there exists $s\in {\rm SHom}_{\tilde{R}_2}(\tilde{R}_1\otimes_{\qc}\tilde{R}_2,\tilde{R}_2)$ with image in ${\rm SDiff}^{\wedge}_{I}(\tilde{R}_1,\tilde{R}_2)$. 

Moreover, if $\tau':\tilde{R}_1\xr{} \tilde{R}_2$ is an isomorphism of $D$-algebras inducing the identity modulo $I$, then we can find an $s$ such that $(s\circ \delta_0)\otimes_D D/(q-1)^{p-1}=\tau'\otimes_{D} D/(q-1)^{p-1}$.
\begin{proof}
Let $m=p-1$, there is $u\in \Z_p\lf q-1 \rf$ such that $u\equiv 1 \mod q-1$ and $u\cdot \dd -p \in (q-1)^m\Z_p\lf q-1 \rf$; we set $d=u\cdot \dd$ and $x=q-1$. Note that $d-p\in (x^m)$ and $\phi(x)\subset (x\cdot d)$. 

Let $J\subset \tilde{R}_1\otimes_{\qc}\tilde{R}_2$ be the $q$-PD ideal. We will define maps $\gamma_{p^k}: J\xr{} J$ by induction on $k$. For $k=0$, we take $\gamma_1={\rm id}_J$. These maps will satisfy $\phi(\gamma_{p^k}(J))\subset d^{m^k} (\tilde{R}_1\otimes_{\qc}\tilde{R}_2)$. 

We define
\begin{align}
\gamma_{p^k}(a) &=  -\delta(\gamma_{p^{k-1}}(a)) + \frac{d^{m^{k-1}} - (d-p)^{m^{k-1}}}{p} \cdot \frac{\phi(\gamma_{p^{k-1}}(a))}{d^{m^{k-1}}} \label{align-gamma-1} \\
                &= \frac{1}{p}\cdot \left( \gamma_{p^{k-1}}(a)^p - (d-p)^{m^{k-1}}  \cdot \frac{\phi(\gamma_{p^{k-1}}(a))}{d^{m^{k-1}}}\right), \label{align-gamma-2}
\end{align}
and have to show that $\phi(\gamma_{p^k}(J))\subset d^{m^k}$ and $\gamma_{p^k}(J)\subset J$. Now, $\phi(\gamma_{p^k}(J))\subset d^{m^k}$ follows from \eqref{align-gamma-2}, $d-p\in (x^m)$, and $(p,d)$ are a regular sequence. For $\gamma_{p^k}(J)\subset J$, we note that by using $\frac{d^{m^{k-1}} - (d-p)^{m^{k-1}}}{p} \equiv d^{m^{k-1}-1} \mod x$ and $x\in J$, we only have to show 
$
-\delta(a) + \frac{\phi(a)}{d}\in J 
$ 
for all $a\in J$. Since $u\equiv 1 \mod q-1$, this follows from the definition of a $q$-PD pair. 

Find coordinates $x_1,\dots, x_n$ for $\tilde{R}_1$ and write $\tilde{R}_1\otimes_{{\rm inf}} \tilde{R}_2 = R_2\lf \epsilon_1,\dots, \epsilon_n \rf$, with  $\delta_0(x_i)=\tau'(x_i)+\epsilon_i$. We can construct a basis $(\Gamma_I)$ for $\tilde{R}_1\otimes_{\qc}\tilde{R}_2$ and a section $s$ as was done in the proof of Proposition \ref{proposition-retractions} and by using the newly defined $\gamma_{p^k}$. 

Evidently, $(s\circ \delta_0)\otimes_D D/x^{m} = \tau'\otimes_{D} D/x^{m}$ holds. 
Moreover, it is not hard to find $n_I\in \Z_p\lf x \rf $ such that $n_I$ is power of $p$ modulo $x$, $n_I\cdot \Gamma_I\in R_2\lf \epsilon_1,\dots, \epsilon_n \rf$, and $n_I\cdot \Gamma_I =  \epsilon^I$ modulo $x^m$, where we used the notation $\epsilon^I=\epsilon_1^{I_1}\cdots \epsilon_n^{I_n}$. 

Writing $n_I\Gamma_I = \sum_{J} B_{J,I}\cdot \epsilon^J$, we get $B_{I,I}=1$, $B_{J,I}\in (x^m)$ for $I\neq J$, and, by using \ref{align-gamma-2}, $\lim_I B_{J,I}=0$ for the $x$-adic topology. Now, writing $\epsilon^I = \sum_{J} A_{J,I}\cdot n_J\Gamma_J$, it is not hard to conclude $\lim_I A_{J,I} =0$ for all $J$. Since $s(\epsilon^I)=A_{0,I}$, we are done.
\end{proof}
\end{proposition}

\subsection{A global approach} \label{subsection-only-subsection?-good-lord-if-I-had-more-time}

In this subsection we are going to patch the functors from Corollary \ref{corollary-equivalance-q-connections} for various primes $p$ together in order to obtain a global result on $q$-connections. By definition, the category of $q$-connections depends on the choice of coordinates. The study of their independence was initiated by Scholze \cite{S}. Our approach will not be able to include the prime $2$, but we have no doubt that other approaches can prove the results in full generality.

\newcommand{\xx}{q-1}
\newcommand{\DD}{\Z\lf \xx\rf}
\newcommand{\Dp}{\Z_p\lf \xx\rf}

Let $R$ be smooth over $\Z$. We set $R_{q}=R\lf q-1 \rf$, and denote by $R_p$ and $R_{q,p}$ the $p$-adic completions. If $\psi: \Z[x_1,\dots, x_n]\xr{} R$ is \'etale, then we denote by $R_{q,p,\psi}$ the $\Z_{p}\lf q-1 \rf$-algebra $R_{q,p}$ together with the $\delta$-structure induced by $\phi(\psi(x_i))=x_i^p$. We note that $(R_{q,p,\psi}, (\xx))$ is an object in $(R_p/\Dp)_{\qc}$.

As in the $p$-adic setup, we set $R_q\otimes_{{\rm inf}} R_q=R_q\lf \epsilon_1,\dots,\epsilon_n\rf$, $\delta_1(a)=a$, for all $a\in R_q$, and $\delta_0$ is induced by $\delta_0(\psi(x_i))=\psi(x_i)+\epsilon_i$ for some $\psi$. Again, we note that $R_q\otimes_{{\rm inf}} R_q$ together with the ideal $J=(q-1,\epsilon_1,\dots, \epsilon_n)$, $\delta_0$, and $\delta_1$, is idependent of $\psi$. It is the coproduct of $R_q$ with itself in a suitable category of infinitisimal thickenings of $R$ over $\DD$.     

We set 
$$
{\rm Diff}^{\wedge}_{(\xx)}(R_q,R_q)=\{ s\in \Hom_{R_q}((R_q\otimes_{{\rm inf}} R_q)_{\delta_1}, R_q)\mid \lim_n s(J^n)=0 \quad \text{$\xx$-adically} \},
$$
and note that 
$$
{\rm Diff}^{\wedge}_{(\xx)}(R_q,R_q) \xr{} \Hom_{\DD}(R_q,R_q), \quad s\mapsto s\circ \delta_0
$$
is injective. 

Depending on $\psi$, we have the differential operators $\nabla_{x_i,q}\in {\rm Diff}^{\wedge}_{(\xx)}(R_q,R_q)$ from \cite{S}. We will use the notations 
$
\nabla^I_{\psi,q}:=\prod_{k=1}^n \nabla^{I_k}_{x_k,q}
$, and 
$$
\mc{A}_{\psi}:= \{ \sum_{I} a_I \cdot \nabla^I_{\psi,q} \mid a_I\in R_q, \lim_I a_I=0 \quad \text{$\xx$-adically}\}.
$$
It is a subalgebra of ${\rm Diff}^{\wedge}_{(\xx)}(R_q,R_q)$. 

We have the $p$-adic completion map 
$$
{\rm Diff}^{\wedge}_{(\xx)}(R_q,R_q) \xr{} {\rm Diff}^{\wedge}_{(\xx)}(R_{q,p},R_{q,p}),
$$
and define $\mc{A}_{\psi,p}$ in a similar way. 

\begin{lemma}\label{lemma-A-psi}
We have 
\begin{align}
&\qdiff(R_{q,p,\psi})= \prod_{I\in \Z^n_{\geq 0 }} R_{q,p}\cdot \nabla^I_{\psi,q}, \label{align-A-1} \\
&\mc{A}_{\psi,p} = \{ s\in \qdiff(R_{q,p,\psi}) \mid s\in {\rm Diff}^{\wedge}_{(\xx)}(R_{q,p},R_{q,p})\}. \label{align-A-2}
\end{align}
\begin{proof}
We know $\dd \cdot x_i^{p-1} \cdot \phi\circ (\nabla_{x_i,q} \circ \delta_0) = (\nabla_{x_i,q} \circ \delta_0) \circ \phi$. This implies 
\begin{equation}\label{equation-nabla-phi}
\nabla_{x_i,q} \circ \phi = \dd \cdot x_i^{p-1} \cdot \phi\circ \nabla_{x_i,q} 
\end{equation}
as equality of maps 
$
R_{q,p} \otimes_{{\rm inf}} R_{q,p} \xr{} R_{q,p}.
$ From \eqref{equation-nabla-phi} we see that $\nabla_{x_i,q}$ factors over $R_{q,p,\psi} \otimes_{\qc} R_{q,p,\psi}$. Since $\nabla_{x_i,q} \equiv \nabla_{x_i} \mod q-1$, we have proved \eqref{align-A-1}. 

In fact, it is not  hard to compute the basis $(\Gamma_I)$ dual to $(\nabla_{\psi,q}^I)$. For $[I]_q!:=\prod_{k=1}^n [I_k]_q!$, we get $[I]_q! \cdot \Gamma_I \in R_{q,p} \otimes_{{\rm inf}} R_{q,p}$ and $[I]_q! \cdot \Gamma_I \in J^{\sum_{k=1}^n I_k}$. This proves \eqref{align-A-2}.    
\end{proof}
\end{lemma}

\begin{lemma}\label{lemma-patching}
The map
$$
{\rm Diff}^{\wedge}_{(\xx)}(R_q,R_q)/\mc{A}_{\psi} \xr{} \bigoplus^{\wedge}_{\text{$p$ prime}} {\rm Diff}^{\wedge}_{(\xx)}(R_{q,p},R_{q,p})/\mc{A}_{\psi,p}, 
$$
where $\bigoplus^{\wedge}$ means the $q-1$-adic completion of the direct sum, is bijective. 
\begin{proof}
Let us consider the differential operators of finite rank: 
\begin{align*}
{\rm Diff}_{(\xx)}(R_q,R_q) &= \{ \sum_{I} a_I \cdot \frac{\nabla_{\psi,q}^I}{[I]_q!} \mid \text{the sum is finite}\}, \\
\mc{A}_{\psi}^{{\rm finite}} &= \{ \sum_{I} a_I \cdot \nabla_{\psi,q}^I \mid \text{the sum is finite}\},
\end{align*}
and similarly for the $p$-adic analogs. Then ${\rm Diff}^{\wedge}_{(\xx)}(R_q,R_q)$ is the $\xx$-adic completion of ${\rm Diff}_{(\xx)}(R_q,R_q)$, and $\mc{A}_{\psi}$ is the $\xx$-adic completion of $\mc{A}_{\psi}^{{\rm finite}}$. And similarly for the $p$-adic analogs. 

It is clear that 
$$
{\rm Diff}_{(\xx)}(R_q,R_q)/\mc{A}^{{\rm finite}}_{\psi} \xr{} \bigoplus_{\text{$p$ prime}} {\rm Diff}_{(\xx)}(R_{q,p},R_{q,p})/\mc{A}^{{\rm finite}}_{\psi,p} 
$$
is bijective. Thus the statement follows  by $q-1$-adic completion after observing that ${\rm Diff}_{(\xx)}(R_q,R_q)/\mc{A}^{{\rm finite}}_{\psi}$ has no $q-1$-torsion.
\end{proof}
\end{lemma}

Following our convention, we define ${\rm SX}=\{s\in {\rm X}\mid s\circ \delta_0 \equiv {\rm id} \mod \xx \}$, where ${\rm X}={\rm Diff}^{\wedge}_{(\xx)}(R_q,R_q)$ or ${\rm X}=\mc{A}_{\psi}$ etc. We note that ${\rm SX}$ is a group via the multiplication on ${\rm X}$.

\begin{corollary}\label{corollary-patching}
The map
$$
{\rm SDiff}^{\wedge}_{(\xx)}(R_q,R_q)/{\rm S}\mc{A}_{\psi} \xr{} \bigoplus^{\wedge}_{\text{$p$ prime}} {\rm SDiff}^{\wedge}_{(\xx)}(R_{q,p},R_{q,p})/{\rm S}\mc{A}_{\psi,p}, 
$$
where $\bigoplus^{\wedge}$ means that each element $(s_p)_p$ has only finitely many components $\neq 1$ after reduction modulo $(q-1)^m$ for every $m$. 
\end{corollary} 

\begin{definition}\label{definition-q-connections-global}
A {\em $q$-connection} on $(R_q,\psi)$ is a $\xx$-derived complete $\DD$-module $M$ together with a morphism of $\DD$-algebras $$\nabla\in \Hom_{\DD}(\mc{A}_{\psi}, \Hom_{\DD}(M,M)).$$  

Morphisms of $q$-connections are morphisms of $\DD$-modules that are compatible with $\nabla$ in the obvious way.
This defines the {\em category of $q$-connections} $\qconn_{(R_q,\psi)}$.
 
Automatically, every $q$-connection is an $R_q$-module and every morphism of $q$-connections is a morphism of $R_q$-modules.
\end{definition}

\begin{proposition}\label{proposition-patching}
Suppose $\frac{1}{2}\in R$. Let $\psi_i:\Z[x_1,\dots,x_n]\xr{} R$ be two \'etale maps. There is a unique $s\in {\rm SDiff}^{\wedge}_{(\xx)}(R_q,R_q)/{\rm S}\mc{A}_{\psi_1}$ such that 
$$
s \cdot {\rm S}\mc{A}_{\psi_1,p} = {\rm SHom}_{R_{q,p,\psi_2}}(R_{q,p,\psi_1}\otimes_{\qc} R_{q,p,\psi_2}, R_{q,p})\cap {\rm SDiff}^{\wedge}_{(q-1)}(R_{q,p},R_{q,p}),
$$
for each prime $p$ (with $R_p\neq 0$).
Moreover, $s \mc{A}_{\psi_1} s^{-1} = \mc{A}_{\psi_2}$. 
\begin{proof}
We know that $${\rm SHom}_{R_{q,p,\psi_2}}(R_{q,p,\psi_1}\otimes_{\qc} R_{q,p,\psi_2}, R_{q,p})\cap {\rm SDiff}^{\wedge}_{(q-1)}(R_{q,p},R_{q,p})\neq \emptyset$$ by Proposition \ref{proposition-q-1-adic-sections} and $\frac{1}{2}\in R$. 

It is clear that ${\rm SHom}_{R_{q,p,\psi_2}}(R_{q,p,\psi_1}\otimes_{\qc} R_{q,p,\psi_2}, R_{q,p})\cap {\rm SDiff}^{\wedge}_{(q-1)}(R_{q,p},R_{q,p})$ is a right ${\rm S}\mc{A}_{\psi_1,p}$-torsor hence defines an element in ${\rm SDiff}^{\wedge}_{(q-1)}(R_{q,p},R_{q,p})/{\rm S}\mc{A}_{\psi_1,p}$, which is also trivial modulo $(\xx)^{p-1}$ by Proposition \ref{proposition-q-1-adic-sections}. We can use Corollary \ref{corollary-patching} to obtain $s$ and its uniqueness. 

Moreover, $s \mc{A}_{\psi_1} s^{-1} = \mc{A}_{\psi_2}$ follows from $s\mc{A}_{\psi_1,p} s^{-1} = \mc{A}_{\psi_2,p}$ for every $p$.
\end{proof}
\end{proposition}

This implies the analog of Proposition \ref{proposition-morphism-qdiff}.

\begin{corollary}\label{corollary-morphism-qdiff-global}
Every $s$ as in Proposition \ref{proposition-patching} induces an isomorphism 
$$
\psi_s: \mc{A}_{\psi_1} \xr{} \mc{A}_{\psi_2}, \quad \zeta \mapsto s\cdot \zeta\cdot s^{-1}. 
$$
For a second choice $s'$ there is a $\xi\in \mc{A}_{\psi_2}$ such that $\psi_{s'}(\zeta)=\xi \cdot \psi_{s}(\zeta)\cdot \xi^{-1}$ for all $\zeta$. 
\end{corollary}

Corollary \ref{corollary-morphism-qdiff-global} implies the analog of Corollary \ref{corollary-equivalance-q-connections}. This yields a natural equivalence between the categories of $q$-connections for different coordinates.

\section{Appendix}

\subsection{Derived complete modules} 

Our reference for derived completion is \cite[\href{https://stacks.math.columbia.edu/tag/091N}{Tag 091N}]{stacks-project}. 

Let $A$ be a commutative ring and $I$ a finitely generated ideal. We denote the category of $I$-derived complete modules by ${\rm Mod}^{\wedge}_A$. It is a weak Serre subcategory of the category of $A$-modules ${\rm Mod}_A$ \cite[\href{https://stacks.math.columbia.edu/tag/091U}{Tag 091U}]{stacks-project}. In particular, it is an abelian category and ${\rm Mod}_A^{\wedge} \subset {\rm Mod}_A$ is fully faithful and exact. Every classically $I$-complete $A$-module is automatically $I$-derived complete \cite[\href{https://stacks.math.columbia.edu/tag/091T}{Tag 091T}]{stacks-project}.  

For $M,N \in {\rm Mod}^{\wedge}_A$ the $A$-module $\Hom_{{\rm Mod}^{\wedge}_A}(M,N)=\Hom_A(M,N)$ is $I$-derived complete \cite[\href{https://stacks.math.columbia.edu/tag/0A6E}{Tag 0A6E}]{stacks-project}. The functor $N\mapsto \Hom_R(M,N)$ has a left adjoint functor $N\mapsto N \widehat{\otimes}_A M$, which is given by $N\mapsto H^0((M\otimes_A N)^{\wedge})$ with $K\mapsto K^{\wedge}$ being the derived completion \cite[\href{https://stacks.math.columbia.edu/tag/091V}{Tag 091V}]{stacks-project}. Indeed, this follows from 
\begin{multline*}
 \mathrm{Hom}(N, {\rm RHom}(M,P)) = {\rm Hom}(N \otimes^{L} M, P) = {\rm Hom}((N \otimes^{L} M)^{\wedge}, P) \\ =  {\rm Hom}(\tau_{\geq 0} (N \otimes^{L} M)^{\wedge}, P) = {\rm Hom}(N \wotimes M, P), 
\end{multline*}
where we used that derived completion behaves like a left derived functor\cite[\href{https://stacks.math.columbia.edu/tag/0AAJ}{Tag 0AAJ}]{stacks-project}. 

The functor $N\mapsto N \widehat{\otimes}_A M$ is right exact. We have $N \widehat{\otimes}_A M = M \widehat{\otimes}_A N$, and $(N \widehat{\otimes}_A M)\widehat{\otimes}_A P = N \widehat{\otimes}_A (M\widehat{\otimes}_A P)$. 

\subsubsection{Descent}\label{appendix-descent}
Let $B$ be a derived $I$-complete $A$-algebra. Then $M\mapsto M \widehat{\otimes}_A B$ defines a functor 
$
{\rm Mod}^{\wedge}_A \xr{} {\rm Mod}^{\wedge}_B
$      
that is left adjoint to the forgetful functor ${\rm Mod}^{\wedge}_B \xr{} {\rm Mod}^{\wedge}_A$.
Note that the derived completion functor commutes with the forgetful functor from $B$-modules to $A$-modules \cite[\href{https://stacks.math.columbia.edu/tag/0924}{Tag 0924}]{stacks-project}.  In particular, a $B$-module $M$ is $IB$-derived complete if and only if it is $I$-derived complete as an $A$-module. In order to see that $M\mapsto M \widehat{\otimes}_A B$ is the left adjoint, we use 
$$
\mathrm{Hom}_B(M\wotimes_A B, N) = \mathrm{Hom}_B(M\otimes_A B, N) = \mathrm{Hom}_A(M,N). 
$$ 

We define descent data in the usual way, that is, descent data correspond to $(M,\theta)$, where $M$ is a $IB$-derived complete $B$-module and $\theta$ is a morphism of $B$-modules
$$
M \xr{} M \wotimes_{A} B, 
$$     
where $M \wotimes_{A} B$ is a $B$-module via the right factor. Since $M \wotimes_{A} B$ is also a $B$-module via the action on $M$, we get an induced morphism 
$$
M \wotimes_{A} B \xr{} M \wotimes_{A} B.
$$
We require that this map is an isomorphism and the cocycle condition is satisfied. The cocycle condition is 
$$
(\theta \wotimes \mathrm{id}_B) \circ \theta = (\mathrm{id}_M \wotimes \delta^1_{0}) \circ \theta,
$$
as equality of morphism $M\xr{} M \wotimes_{A} B \wotimes_{A} B$ with $\delta^1_{0}: B\xr{} B\wotimes B$ induced by $b\mapsto 1\otimes b$.
We will denote the category of derived complete descent data by $DD_{B/A}^{\wedge}$. 

As usual, we have the base extension functor 
\begin{equation}
\label{appendix-base-extension}
{\rm Mod}^{\wedge}_A \xr{} DD_{B/A}^{\wedge},\quad M\mapsto (M\wotimes_A B, {\rm id}_M\wotimes  \delta^1_{0}).
\end{equation}
It has a right adjoint functor given by 
$$
 DD_{B/A}^{\wedge} \xr{} {\rm Mod}^{\wedge}_A,\quad (M,\theta) \mapsto {\rm ker}(\theta - \imath)
$$ 
where $\imath: M \xr{} M\wotimes_{A} B$ is induced by $m\mapsto m\otimes 1$.

\begin{proposition}\label{proposition-section-descent}
If there is a section $s:B\xr{} A$ in the category of $A$-modules then base extension ${\rm Mod}^{\wedge}_A \xr{} DD_{B/A}^{\wedge}, M\mapsto M\wotimes_A B$ is an equivalence of categories. 

A quasi-inverse is given by $M\mapsto {\rm ker}(\theta - \imath)$. 
\begin{proof}
The strategy of the proof is taken from \cite[\href{https://stacks.math.columbia.edu/tag/08WE}{Tag 08WE}]{stacks-project}.  

Let $(M,\theta)$ be a descent datum. We set $f=\theta$, $g_1=\theta\wotimes {\rm id}_B$ and $g_2={\rm id}_M\wotimes \delta_{0}^1$. We claim that 
\begin{equation}  \xymatrix{ M \ar[r]^-f &   M \wotimes_{A} B \ar@<1ex>[r]^-{g_1} \ar@<-1ex>[r]_-{g_2} &   M \wotimes_{A} B \wotimes_{A} B } 
\end{equation} 
is a split equalizer (see \cite[\href{https://stacks.math.columbia.edu/tag/08WH}{Tag 08WH}]{stacks-project} for the definition). Indeed, we can take $h:  M \wotimes_{A} B  \xr{\theta^{-1}}  M \wotimes_{A} B  \xr{} M$, where the last arrow uses the $B$-module structure, and $i: M \wotimes_{A} B \wotimes_{A} B \xr{} M \wotimes_{A} B$ is induced by ${\rm id}_M$ and the multiplication $B \wotimes_{A} B \xr{} B$, in order to split the equalizer. This means 
\begin{equation}   h \circ f = {\rm id}_M, \quad f \circ h = i \circ g_1, \quad i \circ g_2 = {\rm id}_{M \wotimes_{A} B}. 
\end{equation}
Furthermore, we claim that the equalizer 
$$
\xymatrix{ {\rm ker}(\theta - \imath) \ar[r]^-{\subset} &   M  \ar@<1ex>[r]^-{\theta} \ar@<-1ex>[r]_-{\imath} &   M \wotimes_{A} B} 
$$
is split (in the category ${\rm Mod}^{\wedge}_A$). We note that $M\xr{\theta} M \wotimes_{A} B \xr{{\rm id}_M\wotimes s} M$, which we denote by $h'$, factors through ${\rm ker}(\theta - \imath)$. Indeed, we have a commutative diagram 
$$
\xymatrix{ 
    M \wotimes_{A} B \ar[d]_{{\rm id}_M\wotimes s} \ar@<1ex>[r]^-{g_1} \ar@<-1ex>[r]_-{g_2} &   M \wotimes_{A} B \wotimes_{A} B  \ar[d]^{{\rm id}_M\wotimes {\rm id}_B \wotimes s} \\
   M  \ar@<1ex>[r]^-{\theta} \ar@<-1ex>[r]_-{\imath} &   M \wotimes_{A} B.
}
$$ 

Therefore we may take $h'$ as section $M\xr{} {\rm ker}(\theta - \imath)$. We use ${\rm id}_M\wotimes s: M \wotimes_{A} B \xr{} M$ as the second splitting. The required identities are obvious. 

Since split equalizers remain equalizers after application of the functor $\wotimes_{A} B$, the map
$
{\rm ker}(\theta - \imath) \wotimes_{A} B \xr{} M
$   
is an isomorphism. This shows that base extension is essentially surjective.

Now suppose $(M,\theta)$ is the base extension of an $A$-module $M'$. By inspection of $h'$, we conclude that the natural map $M'\xr{} {\rm ker}(\theta - \imath)$ is surjective. It is also injective, because $M'\xr{} M'\wotimes_A B, m'\mapsto m'\wotimes 1,$ has  ${\rm id}_{M'}\wotimes s$ as section. Hence it is an isomorphism. This finishes the proof. 
\end{proof}
\end{proposition}

\begin{remark}
The functor $\wotimes_A B$ is not exact in general, even if $B$ is an $I$-completely flat $A$-algebra. We do not know whether there is descent for a $I$-completely faithfully flat $A$-algebra. 
\end{remark}

\subsubsection{$I$-completely flat modules}\label{appendix-completely-flat}

Recall from \cite{BS} that a complex $M$ of $A$-modules is {\em $I$-completely flat} if for any $I$-torsion $A$-module $N$, the derived tensor product $M\otimes^L_A N$ is concentrated in degree $0$. This implies in particular that $M\otimes^L_A A/I$ is concentrated in degree $0$ and a flat $A$-module.

In this section we will only consider the case where $I$ is generated by two elements. For more general results we refer to \cite{Y}. The next proposition may very well follow from \cite{Y}. 

\begin{proposition}\label{proposition-derived-completion-and-flatness}
Suppose $I$ is generated by two elements $p,d$. And suppose that $d$ is a non-zero divisor of $A$. If $M$ is a $(p,d)$-completely flat complex of $A$-modules then the derived completion $M^{\wedge}$ is $(p,d)$-completely flat. 
\begin{proof}
As a first step, suppose that $A=A/d$ and $I$ is principal and generated by $p$. We claim that if $M$ is a $p$-completely flat complex then $M^{\wedge}$ is $p$-completely flat. 

Let ${\rm cone}(p)$ be the cone of the multiplication by $p$ endomorphism of $A$.  Because its cohomology is $p$-torsion,  $M\otimes^L {\rm cone}(p)$ is derived $p$-complete. Therefore we get
\begin{equation}\label{equation-cone-p}
{\rm cone}(M^{\wedge}\xr{p} M^{\wedge}) \cong M\otimes^L {\rm cone}(p).
\end{equation}
Let $T$ be a $p$-torsion module. The complex ${\rm cone}(p)\otimes^L T$ has only cohomology in degree $=0,-1$. Since $M$ is $p$-completely flat, the same holds for ${\rm cone}(p)\otimes^L T\otimes^L M$. Now \eqref{equation-cone-p} implies $H^i(M^{\wedge}\otimes^L T)$ for all $i\neq 0$, because $T$ is $p$-torsion. 

Let us consider the case $I=(p,d)$ now.  The $(p,d)$-complete flatness of $M$ implies the $p$-complete flatness of $M\otimes^{L}_A A/d$ as complex of $A/d$-modules. We already know that the derived $p$-completion $(M\otimes^{L}_A A/d)^{\wedge}$ is $p$-completely flat as complex of $A/d$-modules. Since $(M\otimes^{L}_A A/d)^{\wedge}$ is $d$-torsion, it is automatically $(p,d)$-complete. Therefore we have an isomorphism 
\begin{equation}\label{equation-cone-d}
{\rm cone}(M^{\wedge}\xr{d} M^{\wedge}) \cong (M\otimes^{L}_A A/d)^{\wedge}.
\end{equation}
Let $T$ be a $(p,d)$-torsion $A$-module. We have 
\begin{multline*}
(M\otimes^{L}_A A/d)^{\wedge} \otimes_A^L T = (M\otimes^{L}_A A/d)^{\wedge} \otimes_A^L A/d \otimes_{A/d}^L T =  \\  \left((M\otimes^{L}_A A/d)^{\wedge}[-1]\oplus (M\otimes^{L}_A A/d)^{\wedge}\right) \otimes_{A/d}^L T, 
\end{multline*}
and this complex has cohomology in degree $0,-1$ only. In view of \eqref{equation-cone-d}, we conclude $H^i(M^{\wedge}\otimes_A^L T)$ for all $i\neq 0$.
\end{proof} 
\end{proposition}

\begin{lemma}\label{lemma-base-change-flatness}
Let $A$ and $B$ be bounded prisms. Let $M$ be a derived $(p,d)$-complete and $(p,d)$-completely flat $A$-module. Then $M\wotimes_A B$ is a $(p,d)$-completely flat $B$-module and it is the classical $(p,d)$-completion of $M\otimes_A B$.
\begin{proof}
The complex of $B$-modules $M\otimes_A^L B$ is $(p,d)$-completely flat, hence $(M\otimes_A^L B)^{\wedge}$ is $(p,d)$-completely flat by Proposition \ref{proposition-derived-completion-and-flatness}. Now, \cite[Lemma~3.7(2)]{BS} can finish the proof.   
\end{proof}
\end{lemma}

The next lemma has essentially the same proof. 

\begin{lemma}\label{lemma-tensor-product-flatness}
Let $A$ be a bounded prism. Let $M,N$ be derived $(p,d)$-complete and $(p,d)$-completely flat $A$-modules. Then $M\wotimes_A N$ is $(p,d)$-completely flat and it is the classical $(p,d)$-completion of $M\otimes_A N$.
\end{lemma}

\begin{lemma}\label{lemma-descent-flatness}
Let $B$ be a derived $(p,d)$-complete and $(p,d)$-completely flat $A$-algebra. Suppose there is a section $B\xr{} A$ in the category of $A$-modules. Let $M$ be a derived $(p,d)$-complete $A$-module. If $M\wotimes_A B$ is $(p,d)$-completely flat as $B$-module then $M$ is $(p,d)$-completely flat.
\begin{proof}
Since $M$ is a direct summand of $M\wotimes_A B$, it suffices to consider $(M\wotimes_A B) \otimes^L_A T = (M\wotimes_A B) \otimes^L_B (B\otimes_A^L T)$ for $(p,d)$-torsion $A$-modules $T$. Now, $B\otimes^L_A T\cong B\otimes_A T$ is a $(p,d)$-torsion $B$-module and we are done. 
\end{proof}  
\end{lemma}

\bibliography{qdr}

\begin{thebibliography}{{Sta}20}

\bibitem[BO78]{BO}
Pierre Berthelot and Arthur Ogus.
\newblock {\em Notes on {C}rystalline {C}ohomology}.
\newblock Princeton University Press and University of Tokyo Press, 1978.

\bibitem[BS]{BS}
Bhargav Bhatt and Peter Scholze.
\newblock Prisms and prismatic cohomology.
\newblock arXiv:1905.08229.

\bibitem[Gro68]{G}
A.~Grothendieck.
\newblock Crystals and the de {R}ham cohomology of schemes.
\newblock In {\em Dix expos\'{e}s sur la cohomologie des sch\'{e}mas}, volume~3
  of {\em Adv. Stud. Pure Math.}, pages 306--358. North-Holland, Amsterdam,
  1968.
\newblock Notes by I. Coates and O. Jussila.

\bibitem[Sch]{S}
Peter Scholze.
\newblock Canonical $q$-deformations in arithmetic geometry.
\newblock arXiv:1606.01796v1.

\bibitem[{Sta}20]{stacks-project}
The {Stacks project authors}.
\newblock The stacks project.
\newblock \url{https://stacks.math.columbia.edu}, 2020.

\bibitem[Yek20]{Y}
Amnon Yekutieli.
\newblock Weak proregularity, derived completion, adic flatness, and prisms,
  2020.
\newblock arXiv:2002.04901.

\end{thebibliography}
\end{document}